\documentclass[12pt]{article}

\usepackage[english]{babel}
\usepackage{times}
\usepackage{amsfonts}
\usepackage{amsthm}
\usepackage{amsmath}
\usepackage{amssymb}
\usepackage{latexsym}
\usepackage{graphicx}
\usepackage{color}
\usepackage{caption}
\usepackage{prodint}
\usepackage[round]{natbib}
\usepackage{bm}
\usepackage{mathabx}
\usepackage{ragged2e}

\newcommand{\rnc}{\renewcommand}
\newcommand{\nc}{\newcommand}
\newcommand{\mrm}{\mathrm}

\nc{\mb}{\mathbb}
\nc{\mac}{\mathcal}
\nc{\E}{\mb{E}}
\nc{\N}{\mb{N}}
\nc{\R}{\mb{R}}
\nc{\Q}{\mb{Q}}
\rnc{\P}{\mrm P}
\rnc{\d}{\mrm d}
\nc{\C}{\mac{C}}
\nc{\D}{\mac{D}}
\nc{\B}{\mac{B}}
\nc{\wh}{\widehat}
\nc{\oPo}{\stackrel{p}{\longrightarrow}}
\nc{\oWo}{\stackrel{w}{\longrightarrow}}
\nc{\oDo}{\stackrel{d}{\longrightarrow}}

\newcommand{\Cov}{cov}

\relax
\allowdisplaybreaks[3]

\newcommand{\nri}{n\rightarrow\infty}

\numberwithin{equation}{section}

\newtheorem{definition}{{\sc Definition}\sc}[section]

\newtheorem{thm}[definition]{Theorem}

\newtheorem{rem}[definition]{Remark}

\newtheorem{cond}[definition]{Condition}

\begin{document}

\linespread{1.3}

    
\title{\Large\bf A Discontinuity Adjustment \\ for Subdistribution Function Confidence Bands \\ Applied to Right-Censored Competing Risks Data}
    
\author{Dennis Dobler \\[1ex] 
}


\maketitle

\vspace{-1cm}

{\centering
{\small Ulm University, \ \
Institute of Statistics, \\
Helmholtzstr. 20, \ \
89081 Ulm, \ \
Germany \\
email: dobler@statistik.tu-dortmund.de \ \ \\
}}

\vspace{0.2cm}

{
\centering

\date{February 3, 2017}

}

\begin{center}
{\bf Summary} \vspace{-.4cm}
\end{center}
\noindent 
The wild bootstrap is the resampling method of choice in survival analytic applications.
Theoretic justifications rely on the assumption of existing intensity functions
which is equivalent to an exclusion of ties among the event times.
However, such ties are omnipresent in practical studies.
It turns out that the wild bootstrap should only be applied in a modified manner
that corrects for altered limit variances and emerging dependencies.
This again ensures the asymptotic exactness of inferential procedures.
An analogous necessity is the use of the Greenwood-type variance estimator for Nelson-Aalen estimators
which is particularly preferred in tied data regimes.
All theoretic arguments are transferred to bootstrapping Aalen-Johansen estimators for cumulative incidence functions in competing risks.
An extensive simulation study as well as an application to real competing risks data of male intensive care unit patients suffering from pneumonia illustrate the practicability of the proposed technique.
  

\vfill
\vfill

\vspace{0cm}


\pagenumbering{gobble}
\thispagestyle{empty}

\newpage

\pagenumbering{arabic}

\allowdisplaybreaks

\section{Introduction}

\label{sec:intro}

In survival analysis, the wild bootstrap is a frequently utilized resampling technique;
cf. \cite{lin97}, \cite{martinussen06}, \cite{beyersmann13}.
Typically, large sample properties are verified by relying on the assumption of existing transition intensities.
In the context of a survival time $T$ this implies the existence of a function
$$\alpha(t) = \lim_{\Delta t \downarrow 0} \frac{1}{\Delta t} P(T \in [t, t + \Delta t] \ | \ T \geq t) $$
or, similarly, the existence of a Lebesgue density of the distribution $P^T$.
Without any doubt, this assumption is too strict in practical applications where events are usually recorded on a discrete time lattice, e.g. on a daily basis.
For estimating unknown limit variances of Nelson-Aalen estimators
it is well-known (e.g. \citealt{abgk93})
that Greenwood-type estimators should be utilized in the presence of ties.
\cite{Allignol_10} even found a general preference (also under left-truncation) for this kind of estimator.

From now on, we assume that there are discrete components in the event time distribution
and that the event, if observed, can be classified to one out of $k$ different causes, i.e. competing risks.
Suppose that $n \in \N$ i.i.d. individuals participate in our study,
but that their observability may be independently right-censored by i.i.d. censoring variables.
Based on these observations, all available collected information is contained in the counting and at risk processes,
\begin{align*}
 N_{ji} (t) & = 1\{ ``\text{ individual $i$ is observed to fail in } [0,t]\text{ due to risk $j$''} \ \} 
 \\
 \text{and} \qquad
 Y_i (t) & = 1\{ ``\text{ individual $i$ is under observation at time } t-\text{''} \ \} 
\end{align*}
respectively, with $j=1,\dots,k$ and $i=1,\dots,n$.
This notation may be used to extend the arguments below to also incorporate independent left-truncation in the sense of \cite{abgk93}, Chapter~III.
Define $H^{uc}_1(t) = E (N_{11}(t))$ and $\bar H(t) = E ( Y_1(t))$.
Now, the sum $N_1(t) = \sum_{i=1}^n N_{1i}(t)$ has the compensator
$$\Lambda_1(t) = \sum_{i=1}^n \Lambda_{1i}(t) = \int_0^t Y(u) \d A_1(u) = \sum_{i=1}^n \int_0^t Y_i(u) \d A_1(u),$$
where $A_1(t) = \int_0^t \bar H(u)^{-1} \d H^{uc}_1(u)$ is the cumulative hazard function for a type 1 event.
Therefore, $(M_1(t) = N_1(t) - \Lambda_1(t))_t$ is a square-integrable martingale.
The Nelson-Aalen estimator for $A_1(t)$ is defined as the counting process integral $\wh A_1(t) = \int_0^t Y^{-1}  \d N_1$, where we let $0/0=0$.
The limit covariance function of the normalized process $W_1 = \sqrt{n} ( \wh A_1 - A_1 )$
is determined by the predictable and optional variation processes of the square-integrable martingales
$M_{1i}: t \mapsto N_{1i}(t) - \Lambda_{1i}(t)$.
In the possibly discontinuous case, these are given by
\begin{align*}
 \langle M_{1i} \rangle (t) & = \Lambda_{1i}(t) - \int_0^t \Delta \Lambda_{1i}(u) \d \Lambda_{1i}(u) \\
 \text{and} \quad [ M_{1i} ] (t) & = N_{1i}(t) - 2 \int_0^t \Delta \Lambda_{1i}(u) \d N_{1i}(u) + \int_0^t \Delta \Lambda_{1i}(u) \d \Lambda_{1i}(u),
\end{align*}
respectively, cf. \cite{abgk93}, Section~II.4.
First, we would like to heuristically discuss again why Greenwood-type variances estimators are a good choice if ties in the event times are possible.
This is as well seen with the help of covariation processes:
Let $A_1^*: t \mapsto \int_0^t J(u) A_1(\d u) / Y(u)$ denote the compensator of the Nelson-Aalen estimator $\wh A_1$, where $J(u) = 1\{ Y(u) > 0\}$.
By the usual martingale theory, we have the optional covariation process
\begin{align*}
 [\sqrt{n}(\wh A_1 - A_1^*)](t) = n \Big[ \int_0^\cdot \frac{J \d M_1}{Y} \Big](t)
  = n \int_0^t \frac{J}{Y^2} \d [M_1]
  = n \sum_{i=1}^n \int_0^t \frac{J}{Y^2} \d [M_{1i}],
\end{align*}
where the last equality follows from the orthogonality of all involved single martingales.
Inserting the above representation for the optional variation process yields
\begin{align*}
 [\sqrt{n} (\wh A_1 - A_1^*)](t) & = n \sum_{i=1}^n \int_0^t \frac{J}{Y^2} \d N_{1i}
  - 2 n \sum_{i=1}^n \int_0^t \frac{J}{Y^2} \Delta \Lambda_{1i} \d N_{1i}
  + n \sum_{i=1}^n \int_0^t \frac{J}{Y^2} \Delta \Lambda_{1i} \d \Lambda_{1i} \\
  & = n \int_0^t \frac{J}{Y^2} \d N_1
  - 2n  \int_0^t \frac{J}{Y^2} \Delta A_1 \d N_1
  + n \int_0^t \frac{J}{Y} \Delta A_1 \d A_1.
\end{align*}
In a last step, we replace the unknown increments of the cumulative hazard function $\Delta A_1(u)$
with Nelson-Aalen increments $\Delta \wh A_1(u) = Y(u)^{-1} \Delta N_1(u)$
and obtain the Greenwood-type variance estimator
\begin{align*}
 \wh \sigma_1^2(t) = n \int_0^t\frac{J}{Y^2} \d N_1 - 2 n \int_0^t \frac{J}{Y^2} \frac{\Delta N_1}{Y} \d N_1
 + n \int_0^t \frac{J}{Y} \frac{\Delta N_1}{Y} \frac{\d N_1}{Y}
 = n \int_0^t \frac{Y - \Delta N_1}{Y^3} \d N_1.
\end{align*}
One could as well proceed similarly by starting with the predictable variation process instead.

Now, for the large sample properties,
let $K > 0$ be a terminal time such that $A_1(K) < \infty$.
A simple application of Glivenko-Cantelli theorems
in combination with the continuous mapping theorem show
that this Greenwood-type estimator is uniformly consistent for $\sigma_1^2(t) = \int_0^t (1 - \Delta A_1(u)) \bar H(u)^{-1} \d A_1(u)$ on $[0,K]$.
In fact, this limit is the asymptotic variance of the normalized Nelson-Aalen estimator
in the general set up where ties are possible; see e.g. \cite{vaart96}, Example~3.9.19, for such a result in the classical survival set up
and Appendix~\ref{sec:app_asy_var} for a detailed calculation for the present competing risks situation.
The subtraction of $\Delta A_1$ in the limit variance is crucially different from the continuous case
in which the limit variance function reduces to $\tilde \sigma^2_1: t \mapsto \int_0^t \bar H^{-1}\d A_1 $.
It is precisely this difference which is the cause for the inconsistency of the usual wild bootstrap procedure applied to Nelson-Aalen estimators and, as a consequence, also when applied to Aalen-Johansen estimators of cumulative incidence functions.

The present article is organized as follows:
In Section~\ref{sec:wbs_nae_univ}, we first discuss the implications of the usual wild bootstrap procedure,
propose a discontinuity adjustment for the wild bootstrap,
and finally present a conditional central limit theorem for this new technique.
In the following Section~\ref{sec:nae_mult_conv}, we show that the Nelson-Aalen estimators for different competing risks are in general asymptotically dependent.
Therefore, we present in a next step an extension of the first proposal for the wild bootstrap adjustment
that guarantees the correct limit dependence structure between the components for different risks.
This technique has some direct implications on resampling the Aalen-Johansen estimator for cumulative incidence functions as these depend on all cause-specific hazard functions and, therefore,
also on their dependencies.
We present conditional central limit theorems corresponding to this set-up in Section~\ref{sec:aje_cif},
where also variance estimators for these Aalen-Johansen estimators
and time-simultaneous confidence bands for cumulative incidence functions are deduced.
The performance of these bands in terms of coverage probabilities is analyzed in a simulation study in Section~\ref{sec:simus}
and there it is compared to the behaviour of confidence bands based on the usual, unadjusted wild bootstrap.
In this connection, we consider different variations of discretization coarseness and discretization probabilities.
All considered resampling techniques are applied to a real data example with competing risks in Section~\ref{sec:data_ex},
where confidence bands for the probability of an alive discharge of male patients with pneumonia from intensive care units are constructed.
We conclude with a small discussion in Section~\ref{sec:disc}.
All proofs and various detailed derivations are presented in  Appendices~\ref{sec:app_asy_var}--\ref{sec:app_asy_var_aje}.

\section{Proposed Wild Bootstrap Adjustment for a Univariate Nelson-Aalen Estimator}
\label{sec:wbs_nae_univ}

The wild bootstrap is typically applied in the following way (see e.g. \citealt{lin97}, \citealt{martinussen06}, \citealt{beyersmann13} or \citealt{bluhmki15} for resampling more general counting process-based Nelson-Aalen estimators):
First, the (independent) martingale increments $\d M_{1i}(u)$ in the martingale representation of $\sqrt{n}( \wh A_1 - A_1)$
are replaced by independently weighted counting process increments, i.e., by $\xi_i \d N_{1i}(u)$.
Here, the wild bootstrap weights $\xi_i$ are i.i.d. with zero mean and variance 1.
\cite{bluhmki15} argued that the resulting wild bootstrapped Nelson-Aalen estimator given by
$$ \wh {\tilde W}_1(t) = \sqrt{n} \sum_{i=1}^n \xi_i \int_0^t \frac{\d N_{1i}}{Y} $$
constitutes a square-integrable martingale with predictable and optional variation processes
\begin{align*}
 t \longmapsto \langle \wh {\tilde W}_1 \rangle (t) & = n \int_0^t \frac{\d N_1}{Y^2} = \wh {\tilde \sigma}^2_1(t) \\
 \text{and} \qquad t \longmapsto [ \wh {\tilde W}_1 ] (t) & = n \sum_{i=1}^n \xi_i^2 \int_0^t \frac{\d N_{1i}}{Y^2} = \wh {\wh {\tilde \sigma}}^2_1(t).
\end{align*}
When utilized as variance estimators, both variation processes in the previous display are uniformly consistent for $\tilde \sigma^2$.
Unfortunately, this is even true in the discontinuous case where $\sigma^2 \neq \tilde \sigma^2$.
Therefore, the wild bootstrap procedure, as it is presently applied in the literature,
fails to be consistent for the correct limit distribution of the Nelson-Aalen estimator
whenever ties are an inevitable phenomenon in the analyzed data set.
This is exactly the reason why this wild bootstrap approach consistently overestimates the true variance in case of ties.
This phenomenon does not vanish as $n \rightarrow \infty$
and it is stronger pronounced for coarser lattices of discrete event times;
see the simulation study in Section~\ref{sec:simus}.
This problem occurs while resampling Aalen-Johansen estimators for cumulative incidence functions in competing risks data 
and it occurs as well while resampling general Nelson-Aalen estimators
if ties are present.

This non-trivial problem calls for a general solution.
In the present article, we exemplify the subsequent solution in the right-censored competing risks set-up.
Extensions and modifications to wild bootstrap versions of more general Nelson-Aalen estimators or of Aalen-Johansen estimators in general Markovian situations 
may be obtained in a similar manner, but the limit variances will be much more complicated.
The crucial defect in the wild bootstrap resampling scheme described above is
that the martingale increments $\d M_{1i}(u)$ should not be replaced by $\xi_i \d N_{1i}(u)$
but rather by something which -- considered again as a martingale -- reproduces the correct covariance structure.
Therefore, we replace
$$ \d M_{1i}(u) \quad \text{by} \quad \xi_i \sqrt{\frac{Y(u) - \Delta N_1(u)}{Y(u)}} \d N_{1i}(u), $$
which results in the following wild bootstrap resampling version of the normalized Nelson-Aalen estimator:
\begin{align*}
 \wh W_1(t) = \sqrt{n} \sum_{i=1}^n \xi_i \int_0^t \sqrt{\frac{Y(u) - \Delta N_1(u)}{Y(u)}} \frac{\d N_{1i}(u)}{Y(u)}.
\end{align*}
In a way similar to \cite{bluhmki15}, one can show that $( \wh W_1(t))_{t \in [0,K]}$ is a martingale with respect to the filtration $(\mac F_t)_{t \in [0,K]}$ given by
$$ \mac F_t = \sigma( \xi_i N_{1i}(u), N_{1i}(v), Y_i(v): u \in [0,t], v \in [0,K], i=1,\dots, n ). $$
Analogously, it is easy to see that its predictable and optional variation processes are given by
\begin{align*}
 \wh \sigma^2_1 : t \longmapsto \langle \wh W_1 \rangle (t) = & n \int_0^t \frac{Y(u) - \Delta N_1(u)}{Y^3(u)} \d N_1(u) \\
 \text{and} \quad \wh {\wh \sigma}^2_1 : t \longmapsto [ \wh W_1 ] (t) = & n \sum_{i=1}^n \xi_i^2 \int_0^t \frac{Y(u) - \Delta N_1(u)}{Y^3(u)} \d N_{1i}(u).
\end{align*}
It is well-known that the Greenwood-type variance estimator $\wh \sigma_1^2$ is uniformly consistent for $\sigma^2_1$.
Assuming the existence of the fourth moments of $\xi_i$, a simple application of Chebyshev's inequality
shows the (conditional) consistency of the second estimator; the uniformity in the conditional convergence in probability follows from a P\'olya-type argument.
The conditional weak convergence of the finite-dimensional marginal distributions
of the wild bootstrapped Nelson-Aalen process
follows easily by an application of Theorem~A.1 in \cite{beyersmann13}.
This also shows that the proposed wild bootstrap approach succeeds in maintaining the correct asymptotic covariance function
which had been our aim in the first place.

Denote by $D_1 = \{ t \in [0,K]: \Delta A_1(t) > 0 \}$ the subset of time points for which ties among the type $1$ events are possible.
Similarly, we define $D_2, \dots, D_k$  for the other risks.
Throughout the article, we assume the following technical condition in order to conclude the conditional tightness of the wild bootstrapped Nelson-Aalen process.
\begin{cond}
\label{cond:main}
 $D = \bigcup_{j=1}^k D_j$ has a finite cardinality.
\end{cond}
In practical applications, this assumption is naturally satisfied:
A finite end-of-study time and measurements on a daily or weekly basis result in a finite lattice.
%
%
%
A proof of conditional tightness finally yields the following conditional central limit theorem for the Nelson-Aalen process:
\begin{thm}
 \label{thm:main}
 Assume Condition~\ref{cond:main}.
 Given $\mac F_0$ and as $\nri$, we have the following conditional  weak convergence
 $$ \wh W_1 \oDo U_1 \sim \text{Gauss}(0, \sigma^2_1) \qquad \text{in outer probability} $$
 on the c\`adl\`ag function space $D[0,K]$ equipped with the supremum distance topology,
 where $U_1$ is a Gaussian zero-mean martingale with variance function $t \mapsto \sigma^2_1(t)$.
 That is, the modified wild bootstrap succeeds in reproducing the same limit process of the Nelson-Aalen process, in particular, if ties are present.
\end{thm}
The proof is given in Appendix~\ref{sec:app_wbs_univ_nae}.
When resampling a functional of a multivariate Nelson-Aalen estimator such as the Aalen-Johansen estimator,
it is mandatory to take the covariance structure between all cause-specific Nelson-Aalen estimators into account.
In order to reflect this in the resampling scheme, a further adjustment needs to be done as shown in the following section.
We conclude this section with an application of the present approach to the Kaplan-Meier estimator.
\begin{rem}
 Consider the case of only one risk, i.e. $k=1$ and $W = W_1$.
 The Kaplan-Meier estimator for the survival function $S(t) = P(T > t)$ is
 $$ \wh S(t) = \prod_{0 \leq u \leq t} (1 - \wh A(\d u)), \quad t \in [0,K]. $$
 It exhibits the martingale representation
 $$W_{S}(t) = \sqrt{n}(\wh S(t) - S(t)) = S(t) \int_0^t \frac{W(\d u)}{1 - \Delta A(u)} + o_p(1).$$
 Thus, the discontinuity-adjusted wild bootstrapped normalized Kaplan-Meier estimator is
 $$ \wh W_S = \wh S \int_0^\cdot \frac{\wh W(\d u)}{1 - \Delta \wh A(u)} $$
 and Theorem~\ref{thm:main} in combination with the continuous mapping theorem yields the correct limit process distribution, i.e. a zero-mean Gaussian process with covariance function given by
 $$ (s,t) \mapsto S(s) S(t) \int_0^t \frac{\d A(u)}{(1-\Delta A(u)) \bar H(u)}; $$
 cf. e.g. Example~3.9.31 in \cite{vaart96}.
\end{rem}

\section{Extension to the Joint Convergence of Multiple Nelson-Aalen Estimators in Competing Risks Models}
\label{sec:nae_mult_conv}

We now extend the above martingale notation to all $k \in \N$ competing risks,
i.e. to $M_{1i}, \dots, M_{ki}$, $i=1,\dots,n$.
For analyses involving two different cumulative cause-specific hazard functions in a competing risks set-up,
it is also important to take the asymptotic covariance structure of both Nelson-Aalen estimators into account.
In the absolutely continuous case, this asymptotic covariance function vanishes as all Nelson-Aalen estimators are asymptotically independent;
cf. \cite{abgk93}, Theorem~IV.1.2.
In the presence of ties, however, the situation is quite different:
Here we have for the martingales $M_{1i}$ and $M_{2i}$ of Section~\ref{sec:intro} that
\begin{align*}
 \langle M_{1i}, M_{2i} \rangle (t) & = - \int_0^t Y_i(u) \Delta A_1(u) \d A_2(u), \\
 [ M_{1i}, M_{2i} ] (t) & = - \int_0^t Y_i(u) \Delta A_1(u) \d M_{2i}(u)
    - \int_0^t Y_i(u) \Delta A_2(u) \d M_{1i}(u) + \int_0^t Y_i(u) \Delta A_1(u) \d A_2(u),
\end{align*}
cf. the derivations in Section~II.4 in \cite{abgk93}.
Define $W_j = \sqrt{n}(\wh A_j - A_j)$, $j=1, \dots, k$.
The above variation processes are strong evidence that the normalized Nelson-Aalen estimators
$W_1, W_2, \dots, W_k$ are not asymptotically independent anymore in the presence of ties.
This is indeed the case:
\begin{thm}
 \label{thm:nae_asy_dep}
 As $\nri$, we have on the product c\`adl\`ag function space $D^k[0,K]$, equipped with the $\max$-$\sup$-norm, that
\begin{align*}
 (W_1, W_2, \dots, W_k) \oDo (U_1,U_2, \dots, U_k),
\end{align*}
where $U_1,U_2, \dots, U_k$ are zero-mean Gaussian martingales with variance functions
$$ t \mapsto \sigma^2_j(t) = \int_0^t \frac{1 - \Delta A_j(u)}{\bar H(u)} \d A_j(u), \quad j=1,2,\dots, k $$
and covariance functions (for $j \neq \ell$)
$$ (s,t) \mapsto cov(U_j(s), U_\ell(t)) = -  \int_0^{s \wedge t} \frac{\Delta A_\ell(u)}{\bar H(u)} \d A_j(u) =: \sigma_{j \ell}(s \wedge t) . $$
\end{thm}
See Appendix~\ref{sec:app_asy_cov_csh} for a derivation of this asymptotic covariance function.
In order to account for this dependence structure in a joint convergence consideration,
the wild bootstrap of the previous section needs to be adjusted once more.
Therefore, let $ \xi_{j \ell i}, j, \ell = 1, \dots, k, \ i=1, \dots, n$, be i.i.d. random variables with
$E(\xi_{111})=0$ and $E(\xi^2_{111})=1$,
which are also independent of the data.
Denote by $N = \sum_{j=1}^k N_j$ the number of all kinds of events.
Then, the single components of the wild bootstrap version of the multivariate Nelson-Aalen estimator
$(W_1, \dots, W_k)$ are given by
\begin{align*}
 \wh W_j(t) & = \sqrt{n} \sum_{i=1}^n \xi_{jji} \int_0^t \sqrt{\frac{Y(u) - \Delta N(u)}{Y(u)}} \frac{\d N_{ji}(u)}{Y(u)} \\
  & + \frac{1}{\sqrt{2}} \sum_{\ell =1}^k sign(\ell - j) \Big[ \sqrt{n} \sum_{i=1}^n \xi_{j \ell i} \int_0^t \sqrt{\frac{\Delta N_j(u)}{Y(u)}} \frac{\d N_{\ell i}(u)}{Y(u)}
  +  \sqrt{n} \sum_{i=1}^n \xi_{\ell j i} \int_0^t \sqrt{\frac{\Delta N_\ell (u)}{Y(u)}} \frac{\d N_{j i}(u)}{Y(u)} \Big],
\end{align*}
where $sign(x) = 1\{ x > 0\} - 1 \{ x < 0 \}$ is the signum function.
This signum function is important in order to insure the required negative covariance between all components.
The following large sample properties hold:
\begin{thm}
 \label{thm:nae_asy_dep_wbs}
 Assume Condition~\ref{cond:main}. Given $\mac F_0$ and as $\nri$, we have the following conditional weak convergence on the product c\`adl\`ag function space $D^k[0,K]$, equipped with the $\max$-$\sup$-norm:
\begin{align*}
 (\wh W_1, \wh W_2, \dots, \wh W_k) \oDo (U_1,U_2, \dots, U_k) \qquad \text{in outer probability},
\end{align*}
where $(U_1,U_2, \dots, U_k)$ is the same Gaussian martingale as in Theorem~\ref{thm:nae_asy_dep}.
\end{thm}
The proof is given in Appendix~\ref{sec:app_wbs_mult_nae}.
Note that, if we are interested in just a single univariate Nelson-Aalen estimator,
the present approach yields the same limit distribution as the wild bootstrap technique proposed in Section~\ref{sec:wbs_nae_univ}.
Hence, it does -- asymptotically -- not matter which of both techniques is applied to the univariate Nelson-Aalen estimator.

Variance and covariance estimators (also for the wild bootstrap versions) are again motivated by the predictable and optional covariation processes of the involved martingales.
The resulting estimators turn out to be same as those obtained by the plug-in method:
\begin{align*}
 \wh \sigma^2_j(t) & = n \int_0^t \frac{Y(u) - \Delta N_j(u)}{Y^3(u)} \d N_j(u), \qquad j=1,\dots,k, \\
 \wh \sigma_{j\ell}(t) & = - n \int_0^t \frac{\Delta N_j(u)}{Y^3(u)} \d N_\ell(u), \qquad j, \ell=1,\dots,k; j \neq \ell,
\end{align*}
are the usual Greenwood-type (co)variance estimators and
\begin{align*}
 \wh {\wh \sigma}^2_j( \cdot ) & = n \sum_{i=1}^n \xi_{jji}^2 \int_0^\cdot \frac{Y(u) - \Delta N_j(u)}{Y^3(u)} \d N_{ji}(u) \\
  & + \frac{n}2 \sum_{\ell \neq j} \Big[ \sum_{i=1}^n \xi_{j \ell i}^2 \int_0^\cdot \frac{\Delta N_j(u)}{Y^3(u)} \d N_{\ell i}(u)
   + \sum_{i=1}^n \xi_{\ell j i}^2 \int_0^\cdot \frac{\Delta N_\ell(u)}{Y^3(u)} \d N_{j i}(u) \Big]
  , \qquad j=1,\dots,k, \\
 \wh {\wh \sigma}_{j\ell}(\cdot) & = - \frac{n}2 \Big[ \sum_{i=1}^n \xi_{j \ell i}^2 \int_0^\cdot \frac{\Delta N_j(u)}{Y^3(u)} \d N_{\ell i}(u)
  + \sum_{i=1}^n \xi_{\ell j i}^2 \int_0^\cdot \frac{\Delta N_\ell(u)}{Y^3(u)} \d N_{j i}(u) \Big] , \qquad j, \ell=1,\dots,k; j \neq \ell,
\end{align*}
are the optional process-type (co)variance estimators motivated from the wild bootstrap martingale properties.
Assume that all $\xi_{111}$ have finite fourth moments.
By applications of Glivenko-Cantelli theorems in combination with the continuous mapping theorem,
the Greenwood-type (co)variance estimators $\wh \sigma_j$ and $\wh \sigma_{j\ell}$ are shown to be uniformly consistent for $\sigma_j^2$ and $\sigma_{j \ell}$, respectively.
For the wild bootstrap-type (co)variance estimators, we can parallel the arguments in the proof of Theorem~\ref{thm:main},
after first assuming the existence of fourth moments $E(\xi_{111}^4) < \infty$:
In points of continuity of all cumulative hazard functions, i.e. on $[0,K] \setminus D$,
Rebolledo's martingale central limit theorem applies and it also implies the uniform consistency of the optional variation process increments.
In points of discontinuity, which are finitely many by assumption, we approximate $\wh {\wh \sigma}^2_j$ by $\wh \sigma^2_j$
and apply the conditional Chebyshev inequality (given $\mac F_0$) in order to show the negligibility of the differences $ \wh {\wh \sigma}^2_j -  {\wh \sigma}^2_j$ in probability.
The last argument can be repeated for the covariance estimators.
A final application of the continuous mapping theorem yields
\begin{align*}
 {\wh \sigma}^2_j, \ \wh {\wh \sigma}^2_j \oPo \sigma_j \quad
 \text{and} \quad {\wh \sigma}_{j\ell}, \ \wh {\wh \sigma}_{j\ell} \oPo \sigma_{j\ell}
 \quad \text{uniformly on $[0,K]$ in (conditional outer) probability},
\end{align*}
for all $j \neq \ell$ as $\nri$.

\section{Extension to Aalen-Johansen Estimators for Cumulative Incidence Functions in Competing Risks}
\label{sec:aje_cif}

Denote the cumulative incidence functions by $F_j(t), j=1,\dots,k$,
which specify the probabilities to die due to cause $j$ during the time interval $[0,t]$.
For ease of presentation, we consider the situation of $k=2$ competing risks
which is achieved by aggregating all but the first risk to be the second competing risk.
The general results are obtained by replacing $U_2$, $A_2$, and $F_2$ in the representations below
by $U - U_1$, $A - A_1$, and $1 - S - F_1$, respectively.
Here, $A$ denotes the all-cause cumulative hazard function.
Utilizing the functional delta-method in combination with the weak convergence results for the Nelson-Aalen estimator,
we as well obtain a weak convergence theorem for the Aalen-Johansen estimator
$ \wh F_1(t) = \int_0^t \wh S(u-) \d \wh A_1(u)  $
for the cumulative incidence function
$ F_1(t) = \int_0^t S(u-) \d A_1(u) $:
\begin{thm}
\label{thm:aje_asy}
 As $\nri$,
\begin{align*}
 W_{F_1} = \sqrt{n} ( \wh F_1 - F_1 ) \oDo U_{F_1}
 {=} \int_0^\cdot \frac{1 - F_2(u-) - F_1( \cdot )}{1 - \Delta A(u)} \d U_1(u)
 + \int_0^\cdot \frac{F_1(u-) - F_1(\cdot)}{1 - \Delta A(u)} \d U_2(u),
\end{align*}
where $U_{F_1}$ is a zero-mean Gaussian process with covariance function
\begin{align*}
 \sigma^2_{F_1}: (s,t) \mapsto & \int_0^{s \wedge t} \frac{(1 - F_2(u) - F_1(s))(1 - F_2(u) - F_1(t))}{\bar H(u)} \frac{1- \Delta A_1(u)}{(1- \Delta A(u))^2} \d A_1(u) \\
    & + \int_0^{s \wedge t} \frac{(F_1(u) - F_1(s))(F_1(u) - F_1(t))}{\bar H(u)} \frac{1- \Delta A_2(u)}{(1- \Delta A(u))^2} \d A_2(u) \\
    & - \int_0^{s \wedge t} \frac{(1 - F_2(u) - F_1(s))(F_1(u) - F_1(t))}{\bar H(u)} \frac{\Delta A_1(u)}{(1- \Delta A(u))^2} \d A_2(u) \\
    & - \int_0^{s \wedge t} \frac{(1 - F_2(u) - F_1(t))(F_1(u) - F_1(s))}{\bar H(u)} \frac{\Delta A_2(u)}{(1- \Delta A(u))^2} \d A_1(u) .
\end{align*}
\end{thm}
For the application of the functional delta-method, note that the Aalen-Johansen estimator in the present competing risks framework
is a combination of the Wilcoxon and the product integral functional applied to the multivariate Nelson-Aalen estimator.
Both of these functionals are Hadamard-differentiable as shown for example in Section~3.9 of \cite{vaart96}.
A derivation of the above asymptotic covariance function is presented in Appendix~\ref{sec:app_asy_var_aje}.

Now, an appropriate wild bootstrap version of $\sqrt{n} ( \wh F_1 - F_1 )$ is given by
\begin{align*}
 \wh W_{F_1}(t) = \int_0^\cdot \frac{1 - \wh F_2(u-) - \wh F_1( \cdot )}{1 - \Delta \wh A(u)} \d \wh W_1(u)
 + \int_0^\cdot \frac{\wh F_1(u-) - \wh F_1(\cdot)}{1 - \Delta \wh A(u)} \d \wh W_2(u),
\end{align*}
where $\wh W_1$ and $\wh W_2$ are again the wild bootstrap versions of the Nelson-Aalen estimator as presented in Section~\ref{sec:nae_mult_conv}.
Using similar martingale arguments as in Appendix~\ref{sec:app_wbs_univ_nae}, 
we obtain the following conditional central limit theorem for the wild bootstrap version of the Aalen-Johansen estimator:
\begin{thm}
 \label{thm:aje_asy_dep_wbs}
 Assume Condition~\ref{cond:main}. Given $\mac F_0$ and as $\nri$, we have the following weak convergence on the c\`adl\`ag function space $D[0,K]$, equipped with the $\sup$-norm:
\begin{align*}
 \wh W_{F_1} \oDo U_{F_1} \qquad \text{in probability},
\end{align*}
where $U_{F_1}$ is the same Gaussian process as in Theorem~\ref{thm:aje_asy}.
\end{thm}

\begin{rem}[The weird bootstrap]
\label{rem:weird_bs}
 Note that the very same proofs may be applied to verify
 that the above conditional central limit theorems hold for the weird bootstrap as well.
 This resampling scheme corresponds to choosing $\xi_{j \ell i} + 1 \sim Bin(Y(X_i),\max(1,Y(X_i))^{-1})$,
 where $X_i$ is the censoring or event time of individual $i$, whichever comes first.
 This is a particular choice of the data-dependent multiplier bootstrap of \cite{dobler15}.
 In their article, heuristic arguments for the second order correctness under absolute continuity of the data have shown
 that centered unit Poisson variates and weird bootstrap multipliers perform favorably in comparison to standard normal wild bootstrap weights.
 In order to also check the preference of either of the first two resampling procedures in the present set-up, where ties are allowed,
 we included the weird bootstrap yielding competing inference methods into the subsequent simulation study.
\end{rem}

Estimators for $\sigma^2_{F_1}$ and its wild bootstrap variant are obtained similarly as such estimators for the Nelson-Aalen (co)variances, i.e. via plug-in:
\begin{align*}
 \wh \sigma^2_{F_1}: (s,t) \mapsto & \int_0^{s \wedge t} \frac{(1 - \wh F_2(u-) - \wh F_1(s))(1 - \wh F_2(u-) - \wh F_1(t))}{(1- \Delta \wh A(u))^2} \d \wh \sigma_1^2(u) \\
    & + \int_0^{s \wedge t} \frac{(\wh F_1(u-) - \wh F_1(s))(\wh F_1(u-) - \wh F_1(t))}{(1- \Delta \wh A(u))^2}  \d \wh \sigma_2^2(u) \\
    & + \int_0^{s \wedge t} \frac{(1 - \wh F_2(u-) - \wh F_1(s))(\wh F_1(u-) - \wh F_1(t))}{(1- \Delta \wh A(u))^2} \d \wh \sigma_{12}(u) \\
    & + \int_0^{s \wedge t} \frac{(1 - \wh F_2(u-) - \wh F_1(t))(\wh F_1(u-) - \wh F_1(s))}{(1- \Delta \wh A(u))^2} \d \wh \sigma_{21}(u) .
\end{align*}
Similarly, $\wh {\wh \sigma}^2_{F_1}$ is obtained by replacing all estimators $\wh \sigma_j$ and $\wh \sigma_{j \ell}$, $j \neq \ell$,
by their wild bootstrap counterparts $\wh {\wh \sigma}_j$ and $\wh {\wh \sigma}_{j \ell}$, respectively.
Their uniform (conditional) consistencies for $\sigma^2_{F_1}$ on $[0,K]$ follow immediately by the uniform consistency of the Nelson-Aalen (co)variance estimators and the continuous mapping theorem.

\begin{rem}[Deduced confidence bands]
\label{rem:conf_bands}
 Following the lines of \cite{beyersmann13}, time-simultaneous confidence bands for $F_1$ can be deduced.
 In particular, let $\phi(s) = \log ( - \log(1-s))$ be a transformation applied to $F_1$ in order to ensure band boundaries between 0 and 1
 and let $g_1(s) = \log(1 - \wh F_1(s)) / \wh \rho(s)$ and $g_2(s) = \log(1 - \wh F_1(s)) / (1 + \wh \rho^2(s))$ be weight functions leading to the usual equal precision and Hall-Wellner bands, respectively; see \cite{abgk93}.
 Here $\wh \rho^2(s) = \wh \sigma_{F_1}^2(s) / (1 - \wh F_1(s))^2$.
 Let $\wh g_1$ and $\wh g_2$ be their wild bootstrap counterparts,
 i.e. the variance estimates $\wh \sigma_{F_1}^2$ are replaced by $\wh {\wh \sigma}_{F_1}^2$.
 The confidence bands for $F_1$ are then derived from the asymptotics of the supremum distance
 $Z_{1j} = \sup_{u \in [t_1,t_2]} | \sqrt{n} g_j(u) (\phi(\wh F_1(u)) - \phi(F_1(u))) | $
 and its wild bootstrap counterpart
 $\wh Z_{1j} = \sup_{u \in [t_1,t_2]} | \wh g_j(u) \phi'(\wh F_1(u)) \wh W_{F_1}(u) |$, where $[t_1,t_2] \subset [0,K]$ and $j \in \{1,2\}$.
 Let $q_{0.95,j}$ be the conditional $95\%$ quantile of $\wh Z_{1j}$ given the data.
 The resulting asymptotic $95\%$ confidence bands are
 $ 1 - (1 - \wh F_1(s) )^{\exp(\pm n^{-1/2} q_{0.95,j} / g_j(s) )}, \ s \in [t_1,t_2], \ j=1,2. $
\end{rem}

\section{Small Sample Behaviour}
\label{sec:simus}

\begin{figure}[ht]
  \vspace{-1cm}
 \includegraphics[width=.5\textwidth, height=.25\textwidth]{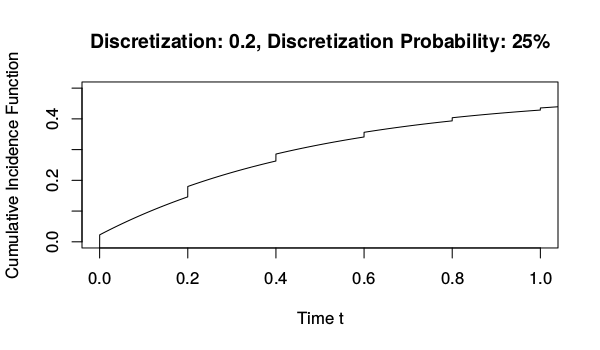}
 \includegraphics[width=.5\textwidth, height=.25\textwidth]{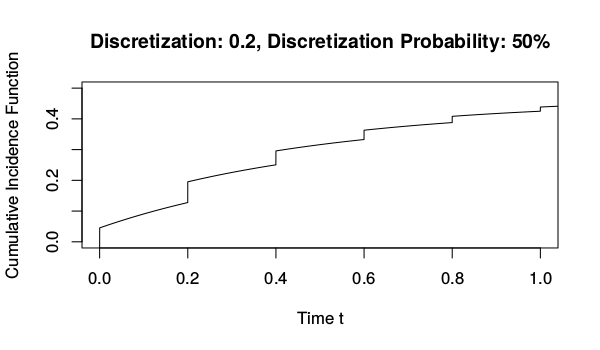} \\[-0.7cm]
 \includegraphics[width=.5\textwidth, height=.25\textwidth]{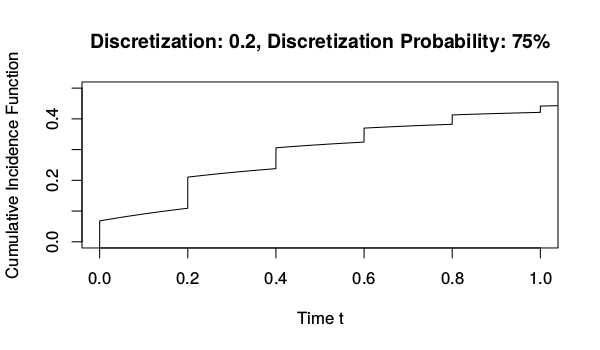}
 \includegraphics[width=.5\textwidth, height=.25\textwidth]{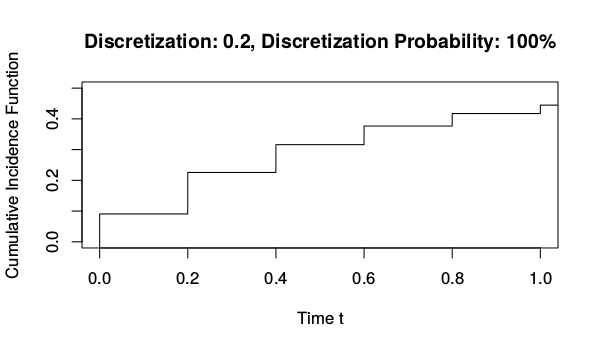} \\[-0.7cm]
 \includegraphics[width=.5\textwidth, height=.25\textwidth]{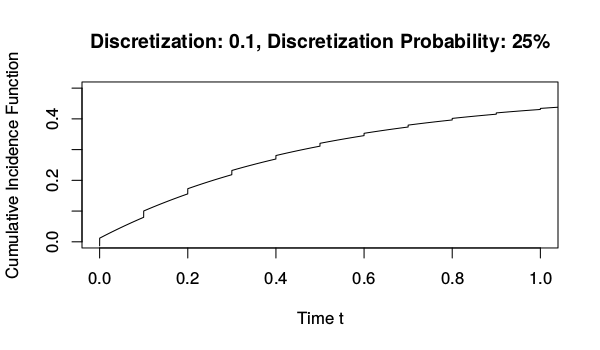}
 \includegraphics[width=.5\textwidth, height=.25\textwidth]{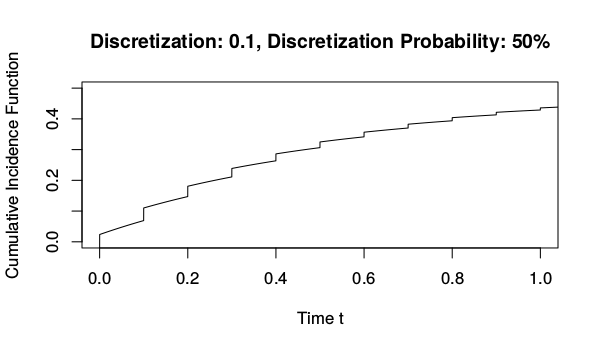} \\[-0.7cm]
 \includegraphics[width=.5\textwidth, height=.25\textwidth]{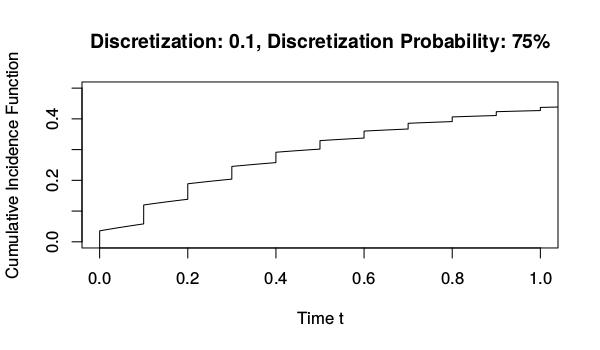}
 \includegraphics[width=.5\textwidth, height=.25\textwidth]{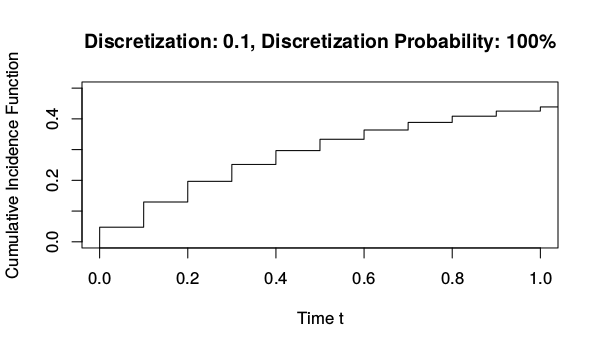} \\[-0.7cm]
 \includegraphics[width=.5\textwidth, height=.25\textwidth]{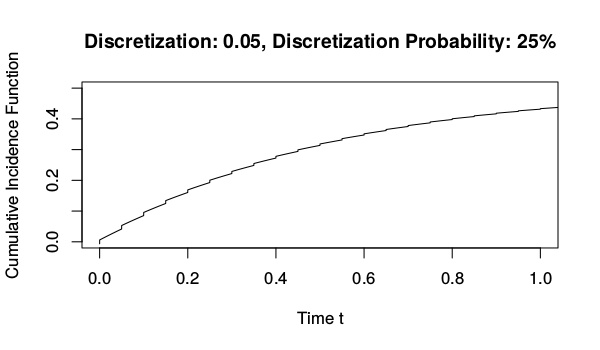}
 \includegraphics[width=.5\textwidth, height=.25\textwidth]{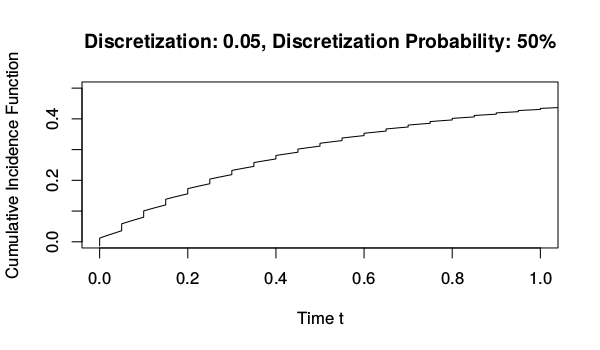} \\[-0.7cm]
 \includegraphics[width=.5\textwidth, height=.25\textwidth]{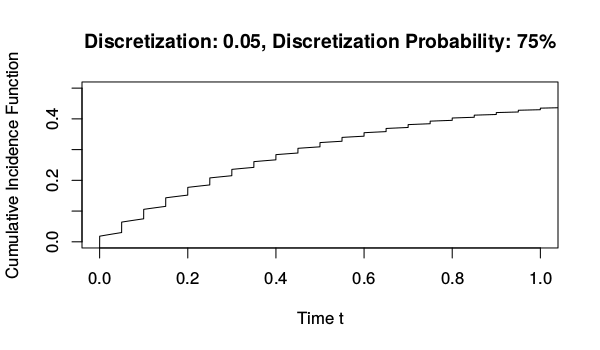}
 \includegraphics[width=.5\textwidth, height=.25\textwidth]{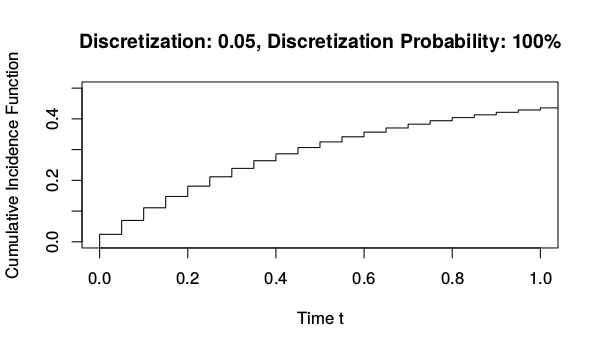}
 \caption{Cumulative incidents functions $F_1^{p,k}$ underlying the present simulations.}
 \label{fig:cifs}
\end{figure}

We empirically assess the difference between the common wild bootstrap approach and the adjusted wild bootstrap proposed in this article via simulation studies.
We simulated the wild bootstrap procedures based on standard normal and centered unit Poisson multipliers as well as the weird bootstrap of Remark~\ref{rem:weird_bs}.
These methods are compared in terms of the simulated coverage probabilities of the confidence bands described in Remark~\ref{rem:conf_bands}.
We consider a simulation set-up motivated by \cite{dobler14},
i.e. we chose the cause-specific hazard rates $\alpha_1(t) = \exp(-t)$ and $\alpha_2(t) = 1 - \exp(-t)$
which yield the cumulative function of the first risk $F_1(t) = 0.5 (1  - \exp(-2t))$.
In order to allow for tied data, we pre-specify different discretization lattices
and round different proportions of the population to the nearest discretization point.
In particular, we choose the discretization lattices to be $\{0, \frac1k, \frac2k, \dots \}$, where $k \in \{5, 10, 20\}$,
and the discretization probabilities to be $p \in \{0, 0.25, 0.5, 0.75, 1\}$.
The resulting theoretic cumulative incidence functions
$$ F_1^{p, k} (t) = p F_1\Big( \frac{[ kt - 0.5 ]}{k} + \frac{0.5}{k} \Big) + (1-p) F_1(t) $$
are presented in Figure~\ref{fig:cifs},
with the exception of the continuous function $F_1$.
Here $[ s ]$ denotes the integer closest to $s \in \R$.
For simulating data, which have the desired cumulative incidence function $F_1^{p,k}$,
it is mandatory to first round the event times $T_i$, and then generate the event types $\varepsilon_i$ in a second step, according to the formula
$$ P \Big( \varepsilon_i = 1 \Big| \frac{[T_i k]}{k} = u \Big) = \frac{F_1(u + \frac{1}{2k}) - F_1(\max(u - \frac{1}{2k},0))}{S(\max(u - \frac{1}{2k},0)) - S( u + \frac{1}{2k})}, $$
where $S: t \mapsto \exp(-t)$ denotes the survival function of the continuous random variables $T_i$.

Censoring is introduced by i.i.d. standard exponentially distributed random variables.
If the $i$th survival time is discretized, then we discretize the $i$th censoring time as well.
Finally, we take the minimum out of each such pair and mark an individual as censored whenever the (discretized) censoring time precedes the (discretized) survival time.
The sample size increases from $n = 50$ to $n = 250$ in steps of $25$.
We choose the time interval, along which asymptotic $95\%$ confidence bands shall be constructed, to be $[0.25, 0.75]$.
The simulations have been conducted using \verb=R= version 3.2.3 (\citealp{Rteam})
using 10,000 outer Monte Carlo iterations and 999 wild bootstrap replicates.

Tables~\ref{tab:EP5} to~\ref{tab:HW20} contain the simulated coverage probabilities
of equal precision and Hall-Wellner bands for simulation set-ups with a discrete component in the cumulative incidence function, i.e. $p > 0$.
The columns of simulation results corresponding to the common wild / weird bootstrap procedures are entitled \emph{old},
whereas the columns showing the results of the respective adjusted wild / weird bootstrap are entitled \emph{new}.

At first, we start with a discussion on the choice of multipliers.
For equal precision bands and in almost all set-ups, there is a pronounced superiority of the wild bootstrap with centered unit Poisson multipliers and the weird bootstrap over the respective coverage probabilities of the bands based on standard normal weights.
This is true for the common resampling procedures as well as for the proposed adjusted bootstraps.
For Hall-Wellner bands, this superiority is not as much pronounced and sometimes even the confidence bands based on standard normal multipliers yield the most accurate coverage probabilities.
But in cases, where this is so, the deviance is only very small.
The phenomenon, that standard normal multipliers yield a worse performance than those with skewness equal to one,
is in line with the findings in a revised version of \cite{dobler15}
where also heuristic theoretic arguments for a second-order correctness of both superior resampling procedures are provided.
As there is, all in all, not much of a difference between the simulated coverage probabilities of the centered unit Poisson wild bootstrap and the weird bootstrap,
we only focus on the results of the Poisson choice.
In general, the equal precision bands are more accurate than the Hall-Wellner bands.

Here, the discretization-adjusted wild bootstrap almost always yielded coverage probabilities closer to the nominal level in comparison to the unadjusted wild bootstrap.
The deviances between these coverage probabilities of each of those two resampling procedures
appears to be larger the higher the discretization probability $p$ and the coarser the discretization lattice is.
For instance, this difference even amounts to $4.03$ percentage points in case of the Hall-Wellner bands, $n=50$, and $p=1$
and to $3.83$ percentage points in case of equal precision bands and the same $p$ and $n$.

In case of $k \in \{5,10\}$, the coverage probabilities of the common wild bootstrap do not appear to converge at all towards $95\%$ as the sample size increases.
Instead, the simulated probabilities fluctuate around $93\%$ or even $92\%$.
On the other hand, the discretization-adjusted equal precision wild bootstrap bands yield much better coverage probabilities
which are greater than $94\%$ or at least in the high $93\%$-region for larger sample sizes.
In contrast to the unadjusted procedure,
we observe for small samples and for the adjusted confidence bands coverage probabilities closer to the nominal level for higher discretization probabilities $p$.
This is only reasonable, as $p=100\%$ corresponds to a multivariate, but not an infinite-dimensional statistical problem.
We do not see this tendency for the unadjusted procedure in case of $k \in \{5, 10\}$,
which again stresses that it is not suitable for these kinds of tied data regimes.

Finally, Table~\ref{tab:HWEPcont} shows the corresponding results for the scenario in which the continuous $F_1$ is the true cumulative incidence function and the usual wild bootstrap technique yields asymptotically exact inference procedures.
However, it is surprising to see that even here the adjusted wild bootstrap yields more accurate confidence bands in comparison to the unadjusted procedure.
Therefore, there is apparently no loss at all in utilizing the discretization adjustment -- with or without ties in the data.

All in all, we conclude that the proposed discontinuity adjustment should always be applied
in order to greatly improve the coverage probabilities of confidence bands for $F_1^{p,k}$.
The present simulation results show this improvement, which amounts to up to two or three percentage points for smaller samples,
in many conducted simulation scenarios.
As the standard normal variate-based wild bootstrap disappoints in general, our final advice is to combine the present discontinuity adjustment with the wild bootstrap based on the Poisson-distributed random variables or the weird bootstrap.
Additionally, equal precision bands should be preferred to Hall-Wellner confidence bands
due to the slight but frequent difference in coverage probabilities.

\begin{table}[ht]
\centering
\scriptsize
\begin{tabular}{c|c|cc|cc|cc}
  \multicolumn{2}{c|}{set-ups} & \multicolumn{2}{c|}{$N(0,1)$} & \multicolumn{2}{c|}{$Poi(1)-1$}  & \multicolumn{2}{c}{weird}  \\ \hline
  $p$ & $n$ & old & new & old & new & old & new \\ \hline
   & 50 & 87.31 & 88.47 & 89.44 & 90.61 & 88.93 & 90.36 \\
   & 75 & 89.48 & 90.39 & 91.55 & 92.24 & 91.13 & 92.01 \\
   & 100 & 90.64 & 91.26 & 92.32 & 92.96 & 92.21 & 92.64 \\
   & 125 & 91.21 & 91.96 & 92.36 & 93.07 & 93.43 & 93.76 \\
  25\% & 150 & 91.12 & 91.66 & 92.48 & 92.98 & 93.37 & 93.75 \\
   & 175 & 91.69 & 92.14 & 93.10 & 93.21 & 93.20 & 93.63 \\
   & 200 & 91.99 & 92.46 & 93.31 & 93.80 & 93.63 & 94.05 \\
   & 225 & 92.51 & 92.98 & 93.29 & 93.67 & 93.58 & 94.03 \\
   & 250 & 92.48 & 92.88 & 93.65 & 93.91 & 93.41 & 93.84 \\ \hline
   & 50 & 87.41 & 89.31 & 89.09 & 91.28 & 89.76 & 91.29 \\
   & 75 & 89.35 & 90.99 & 91.28 & 92.71 & 91.03 & 92.59 \\
   & 100 & 90.01 & 91.53 & 92.39 & 93.64 & 91.88 & 93.29 \\
   & 125 & 90.56 & 91.92 & 92.95 & 93.90 & 92.29 & 93.43 \\
  50\% & 150 & 90.59 & 91.93 & 92.44 & 93.38 & 93.04 & 94.05 \\
   & 175 & 91.21 & 92.22 & 92.58 & 93.49 & 92.66 & 93.73 \\
   & 200 & 91.22 & 92.28 & 92.30 & 93.31 & 93.16 & 94.19 \\
   & 225 & 91.44 & 92.61 & 92.70 & 93.73 & 93.25 & 94.16 \\
   & 250 & 91.72 & 92.82 & 93.13 & 94.27 & 93.07 & 94.05 \\ \hline
   & 50 & 87.23 & 90.60 & 89.12 & 91.56 & 89.47 & 92.11 \\
   & 75 & 90.09 & 92.56 & 90.96 & 93.22 & 91.29 & 93.77 \\
   & 100 & 90.58 & 92.86 & 92.13 & 94.23 & 91.59 & 93.83 \\
   & 125 & 90.83 & 93.03 & 92.13 & 94.26 & 91.75 & 94.08 \\
  75\% & 150 & 90.84 & 93.11 & 91.96 & 93.94 & 92.29 & 94.35 \\
   & 175 & 90.83 & 93.13 & 91.66 & 93.87 & 92.25 & 94.03 \\
   & 200 & 90.91 & 93.11 & 92.45 & 94.46 & 92.50 & 94.49 \\
   & 225 & 91.15 & 93.43 & 92.52 & 94.49 & 92.22 & 94.12 \\
   & 250 & 91.76 & 93.76 & 92.38 & 94.36 & 92.73 & 94.53 \\ \hline
   & 50 & 87.71 & 91.98 & 89.08 & 92.91 & 89.56 & 93.19 \\
   & 75 & 88.79 & 92.70 & 90.26 & 93.74 & 90.21 & 93.64 \\
   & 100 & 89.89 & 93.25 & 90.90 & 94.06 & 90.47 & 93.84 \\
   & 125 & 89.76 & 93.14 & 90.92 & 93.88 & 90.44 & 93.64 \\
  100\% & 150 & 89.64 & 93.21 & 90.84 & 93.97 & 91.00 & 93.87 \\
   & 175 & 90.10 & 93.24 & 90.52 & 93.52 & 91.25 & 93.99 \\
   & 200 & 90.27 & 93.55 & 91.06 & 94.00 & 91.09 & 94.12 \\
   & 225 & 90.01 & 93.07 & 91.30 & 94.19 & 91.08 & 94.26 \\
   & 250 & 89.71 & 93.18 & 90.92 & 93.89 & 90.76 & 93.86 \\
   \hline
\end{tabular}
\caption{Simulated coverage probabilities of equal precision bands in per cent where $k=5$}
\label{tab:EP5}
\end{table}

\begin{table}[ht]
\centering
\scriptsize
\begin{tabular}{c|c|cc|cc|cc}
  \multicolumn{2}{c|}{set-ups} & \multicolumn{2}{c|}{$N(0,1)$} & \multicolumn{2}{c|}{$Poi(1)-1$}  & \multicolumn{2}{c}{weird}  \\ \hline
  $p$ & $n$ & old & new & old & new & old & new \\ \hline
   & 50 & 87.72 & 88.64 & 89.37 & 90.69 & 89.87 & 90.84 \\
   & 75 & 89.22 & 90.16 & 91.86 & 92.55 & 91.75 & 92.35 \\
   & 100 & 90.68 & 91.22 & 92.43 & 93.07 & 93.06 & 93.31 \\
   & 125 & 91.13 & 91.61 & 93.05 & 93.50 & 93.42 & 93.79 \\
  25\% & 150 & 91.63 & 92.16 & 93.5 & 93.77 & 93.48 & 93.79 \\
   & 175 & 91.79 & 92.06 & 93.87 & 94.19 & 93.75 & 93.93 \\
   & 200 & 91.93 & 92.37 & 93.84 & 93.94 & 93.70 & 93.88 \\
   & 225 & 92.40 & 92.77 & 93.54 & 93.89 & 93.73 & 93.97 \\
   & 250 & 92.17 & 92.61 & 94.03 & 94.26 & 94.24 & 94.56 \\ \hline
   & 50 & 87.15 & 88.59 & 89.20 & 90.73 & 89.18 & 90.50 \\
   & 75 & 89.78 & 90.80 & 91.46 & 92.47 & 91.18 & 92.22 \\
   & 100 & 90.24 & 91.32 & 92.58 & 93.25 & 92.53 & 93.47 \\
   & 125 & 91.17 & 92.11 & 92.99 & 93.70 & 93.06 & 93.76 \\
  50\% & 150 & 91.63 & 92.52 & 93.25 & 93.98 & 93.17 & 93.95 \\
   & 175 & 92.65 & 93.29 & 92.67 & 93.34 & 93.06 & 93.71 \\
   & 200 & 92.21 & 92.87 & 93.36 & 94.03 & 93.51 & 94.23 \\
   & 225 & 92.20 & 92.98 & 94.07 & 94.80 & 93.09 & 93.96 \\
   & 250 & 92.96 & 93.52 & 93.71 & 94.23 & 93.37 & 93.93 \\ \hline
   & 50 & 87.63 & 89.82 & 89.19 & 91.41 & 88.95 & 90.91 \\
   & 75 & 89.52 & 91.46 & 91.97 & 93.42 & 91.36 & 92.96 \\
   & 100 & 90.30 & 91.90 & 92.45 & 93.76 & 92.52 & 93.95 \\
   & 125 & 91.44 & 92.88 & 92.93 & 94.03 & 92.83 & 94.14 \\
  75\% & 150 & 91.49 & 92.76 & 93.20 & 94.50 & 93.22 & 94.23 \\
   & 175 & 91.50 & 92.76 & 92.98 & 94.18 & 93.66 & 94.63 \\
   & 200 & 92.00 & 93.08 & 93.92 & 94.94 & 93.75 & 94.85 \\
   & 225 & 92.20 & 93.12 & 93.57 & 94.56 & 93.17 & 94.13 \\
   & 250 & 92.49 & 93.78 & 93.51 & 94.55 & 93.69 & 94.81 \\ \hline
   & 50 & 88.28 & 91.10 & 90.21 & 92.75 & 90.31 & 92.68 \\
   & 75 & 89.81 & 92.2 & 91.70 & 93.90 & 91.77 & 93.57 \\
   & 100 & 90.54 & 92.66 & 92.10 & 93.98 & 92.09 & 93.96 \\
   & 125 & 91.12 & 93.19 & 92.58 & 94.36 & 92.74 & 94.47 \\
  100\% & 150 & 91.41 & 93.44 & 92.83 & 94.61 & 92.77 & 94.77 \\
   & 175 & 91.60 & 93.59 & 92.84 & 94.57 & 92.84 & 94.79 \\
   & 200 & 91.49 & 93.46 & 93.03 & 94.76 & 93.25 & 94.97 \\
   & 225 & 92.32 & 93.98 & 92.92 & 94.70 & 93.38 & 94.95 \\
   & 250 & 91.81 & 93.81 & 92.79 & 94.58 & 92.98 & 94.73 \\
   \hline
\end{tabular}
\caption{Simulated coverage probabilities of equal precision bands in per cent where $k=10$}
\label{tab:EP10}
\end{table}

\begin{table}[ht]
\centering
\scriptsize
\begin{tabular}{c|c|cc|cc|cc}
  \multicolumn{2}{c|}{set-ups} & \multicolumn{2}{c|}{$N(0,1)$} & \multicolumn{2}{c|}{$Poi(1)-1$}  & \multicolumn{2}{c}{weird}  \\ \hline
  $p$ & $n$ & old & new & old & new & old & new \\ \hline
   & 50 & 87.21 & 88.03 & 89.93 & 90.94 & 90.15 & 90.69 \\
   & 75 & 89.75 & 90.33 & 92.18 & 92.50 & 91.78 & 92.38 \\
   & 100 & 91.01 & 91.61 & 93.12 & 93.52 & 92.71 & 93.08 \\
   & 125 & 92.09 & 92.29 & 93.99 & 94.18 & 93.34 & 93.67 \\
  25\% & 150 & 91.87 & 92.39 & 93.63 & 93.78 & 93.69 & 93.85 \\
   & 175 & 92.44 & 92.71 & 94.03 & 94.40 & 94.19 & 94.35 \\
   & 200 & 92.76 & 93.03 & 94.18 & 94.35 & 93.76 & 94.07 \\
   & 225 & 93.23 & 93.44 & 94.46 & 94.77 & 94.02 & 94.23 \\
   & 250 & 92.84 & 93.11 & 93.88 & 94.10 & 94.16 & 94.41 \\ \hline
   & 50 & 88.02 & 89.22 & 89.57 & 90.86 & 89.59 & 90.81 \\
   & 75 & 89.88 & 90.74 & 92.37 & 93.28 & 92.34 & 93.00 \\
   & 100 & 90.95 & 91.69 & 93.19 & 93.74 & 93.02 & 93.71 \\
   & 125 & 91.98 & 92.46 & 93.44 & 93.83 & 93.91 & 94.2 \\
  50\% & 150 & 91.90 & 92.54 & 93.60 & 94.01 & 93.97 & 94.27 \\
   & 175 & 92.95 & 93.40 & 94.00 & 94.42 & 94.11 & 94.45 \\
   & 200 & 92.79 & 93.11 & 94.05 & 94.39 & 93.68 & 94.00 \\
   & 225 & 92.78 & 93.25 & 94.16 & 94.53 & 94.28 & 94.67 \\
   & 250 & 93.01 & 93.28 & 93.92 & 94.30 & 93.77 & 94.03 \\ \hline
   & 50 & 88.42 & 89.91 & 90.81 & 92.12 & 90.51 & 91.96 \\
   & 75 & 90.28 & 91.50 & 92.73 & 93.94 & 92.85 & 93.67 \\
   & 100 & 91.43 & 92.45 & 93.30 & 94.21 & 93.02 & 93.85 \\
   & 125 & 91.41 & 92.30 & 93.84 & 94.58 & 93.57 & 94.08 \\
  75\% & 150 & 92.73 & 93.31 & 94.01 & 94.77 & 93.54 & 94.25 \\
   & 175 & 92.44 & 93.14 & 94.44 & 95.04 & 94.09 & 94.82 \\
   & 200 & 93.04 & 93.84 & 94.13 & 94.65 & 94.02 & 94.58 \\
   & 225 & 93.32 & 93.94 & 94.41 & 94.95 & 94.08 & 94.69 \\
   & 250 & 93.03 & 93.59 & 94.27 & 94.85 & 94.29 & 94.77 \\ \hline
   & 50 & 88.87 & 90.73 & 91.04 & 92.79 & 91.16 & 92.69 \\
   & 75 & 90.95 & 92.58 & 93.21 & 94.28 & 93.02 & 94.24 \\
   & 100 & 91.98 & 93.28 & 93.74 & 95.01 & 93.74 & 94.67 \\
   & 125 & 92.91 & 93.94 & 93.69 & 94.88 & 94.13 & 95.13 \\
  100\% & 150 & 92.82 & 93.88 & 94.18 & 95.10 & 94.58 & 95.49 \\
   & 175 & 93.16 & 94.25 & 94.88 & 95.59 & 94.85 & 95.74 \\
   & 200 & 93.53 & 94.51 & 94.38 & 95.30 & 94.59 & 95.29 \\
   & 225 & 94.06 & 94.76 & 94.63 & 95.59 & 94.69 & 95.39 \\
   & 250 & 93.12 & 94.31 & 94.95 & 95.79 & 94.79 & 95.60 \\
   \hline
\end{tabular}
\caption{Simulated coverage probabilities of equal precision bands in per cent where $k=20$}
\label{tab:EP20}
\end{table}

\begin{table}[ht]
\centering
\scriptsize
\begin{tabular}{c|c|cc|cc|cc}
  \multicolumn{2}{c|}{set-ups} & \multicolumn{2}{c|}{$N(0,1)$} & \multicolumn{2}{c|}{$Poi(1)-1$}  & \multicolumn{2}{c}{weird}  \\ \hline
  $p$ & $n$ & old & new & old & new & old & new \\ \hline
   & 50 & 88.65 & 89.93 & 89.33 & 90.79 & 88.85 & 90.18 \\
   & 75 & 89.82 & 90.60 & 90.59 & 91.54 & 90.39 & 91.32 \\
   & 100 & 90.51 & 91.11 & 90.93 & 91.82 & 91.13 & 91.88 \\
   & 125 & 90.92 & 91.69 & 91.11 & 91.82 & 91.07 & 91.92 \\
  25\% & 150 & 91.41 & 91.70 & 91.53 & 92.20 & 92.21 & 92.74 \\
   & 175 & 91.79 & 92.17 & 92.75 & 93.17 & 92.65 & 93.16 \\
   & 200 & 92.24 & 92.67 & 92.18 & 92.59 & 92.26 & 92.73 \\
   & 225 & 92.22 & 92.70 & 92.75 & 93.20 & 92.64 & 93.14 \\
   & 250 & 92.38 & 93.03 & 92.53 & 93.05 & 92.44 & 92.68 \\ \hline
   & 50 & 87.72 & 90.36 & 88.52 & 90.87 & 88.77 & 91.14 \\
   & 75 & 89.17 & 91.06 & 89.54 & 91.43 & 89.92 & 91.87 \\
   & 100 & 89.86 & 91.62 & 90.26 & 91.87 & 90.60 & 92.18 \\
   & 125 & 90.86 & 92.28 & 91.27 & 92.69 & 91.00 & 92.39 \\
  50\% & 150 & 90.99 & 92.28 & 91.34 & 92.77 & 91.15 & 92.57 \\
   & 175 & 91.41 & 92.67 & 91.76 & 92.84 & 91.12 & 92.42 \\
   & 200 & 91.77 & 92.99 & 92.20 & 93.57 & 91.78 & 92.86 \\
   & 225 & 91.70 & 92.87 & 91.72 & 92.95 & 91.94 & 93.14 \\
   & 250 & 91.44 & 92.84 & 92.00 & 93.20 & 91.65 & 93.04 \\ \hline
   & 50 & 87.35 & 90.83 & 88.74 & 92.14 & 88.32 & 91.82 \\
   & 75 & 89.97 & 92.41 & 90.57 & 93.39 & 90.09 & 93.16 \\
   & 100 & 89.89 & 92.60 & 90.32 & 93.29 & 90.19 & 93.01 \\
   & 125 & 90.52 & 92.98 & 90.95 & 93.38 & 91.03 & 93.55 \\
  75\% & 150 & 90.61 & 93.03 & 91.15 & 93.80 & 91.16 & 93.67 \\
   & 175 & 91.03 & 93.41 & 91.54 & 93.95 & 90.88 & 93.37 \\
   & 200 & 91.29 & 93.53 & 90.98 & 93.49 & 91.30 & 93.41 \\
   & 225 & 91.63 & 93.76 & 91.70 & 93.52 & 91.11 & 93.61 \\
   & 250 & 91.44 & 93.55 & 91.55 & 93.78 & 91.95 & 94.24 \\ \hline
   & 50 & 88.99 & 92.88 & 89.60 & 93.67 & 89.68 & 93.41 \\
   & 75 & 90.32 & 93.74 & 90.71 & 94.31 & 90.23 & 93.99 \\
   & 100 & 90.59 & 94.04 & 91.17 & 94.37 & 90.60 & 94.14 \\
   & 125 & 91.27 & 94.59 & 91.79 & 94.81 & 91.56 & 94.55 \\
  100\% & 150 & 91.63 & 94.59 & 91.42 & 94.50 & 91.20 & 94.60 \\
   & 175 & 91.53 & 94.56 & 91.95 & 94.80 & 90.97 & 94.22 \\
   & 200 & 91.78 & 94.65 & 91.82 & 94.68 & 91.76 & 94.68 \\
   & 225 & 91.94 & 94.92 & 91.87 & 94.81 & 91.59 & 94.56 \\
   & 250 & 92.40 & 94.92 & 91.86 & 94.76 & 91.35 & 94.46 \\
   \hline
\end{tabular}
\caption{Simulated coverage probabilities of Hall-Wellner bands in per cent where $k=5$}
\label{tab:HW5}
\end{table}

\begin{table}[ht]
\centering
\scriptsize
\begin{tabular}{c|c|cc|cc|cc}
  \multicolumn{2}{c|}{set-ups} & \multicolumn{2}{c|}{$N(0,1)$} & \multicolumn{2}{c|}{$Poi(1)-1$}  & \multicolumn{2}{c}{weird}  \\ \hline
  $p$ & $n$ & old & new & old & new & old & new \\ \hline
   & 50 & 87.21 & 88.70 & 89.27 & 90.49 & 89.00 & 90.36 \\
   & 75 & 89.23 & 90.08 & 90.91 & 91.57 & 90.65 & 91.51 \\
   & 100 & 90.54 & 90.86 & 91.86 & 92.55 & 91.40 & 91.89 \\
   & 125 & 91.49 & 91.96 & 92.45 & 92.86 & 91.59 & 92.25 \\
  25\% & 150 & 92.52 & 92.80 & 92.31 & 92.94 & 92.18 & 92.65 \\
   & 175 & 92.61 & 92.88 & 92.86 & 93.06 & 92.71 & 93.23 \\
   & 200 & 92.25 & 92.76 & 92.66 & 92.88 & 92.80 & 93.20 \\
   & 225 & 92.92 & 93.41 & 93.02 & 93.32 & 93.00 & 93.38 \\
   & 250 & 92.37 & 92.78 & 92.90 & 93.29 & 93.19 & 93.46 \\ \hline
   & 50 & 87.32 & 89.05 & 88.85 & 90.49 & 89.17 & 90.77 \\
   & 75 & 89.94 & 91.20 & 91.08 & 92.27 & 90.65 & 91.81 \\
   & 100 & 91.24 & 92.37 & 91.08 & 92.32 & 91.51 & 92.59 \\
   & 125 & 91.05 & 92.16 & 91.61 & 92.49 & 91.88 & 92.77 \\
  50\% & 150 & 91.99 & 92.82 & 91.92 & 92.79 & 92.00 & 93.01 \\
   & 175 & 92.37 & 93.10 & 92.53 & 93.21 & 92.66 & 93.38 \\
   & 200 & 92.24 & 92.79 & 92.42 & 93.37 & 92.59 & 93.25 \\
   & 225 & 92.29 & 93.22 & 92.75 & 93.44 & 92.87 & 93.64 \\
   & 250 & 92.05 & 92.73 & 92.70 & 93.38 & 92.05 & 92.74 \\ \hline
   & 50 & 88.15 & 90.60 & 88.32 & 91.01 & 89.49 & 91.71 \\
   & 75 & 90.41 & 92.20 & 91.34 & 92.98 & 90.39 & 92.28 \\
   & 100 & 91.25 & 92.74 & 90.97 & 92.66 & 91.19 & 93.00 \\
   & 125 & 91.76 & 93.21 & 91.99 & 93.34 & 91.57 & 93.09 \\
  75\% & 150 & 91.67 & 92.85 & 91.73 & 93.04 & 91.90 & 93.31 \\
   & 175 & 92.30 & 93.57 & 92.50 & 93.75 & 92.08 & 93.34 \\
   & 200 & 92.25 & 93.46 & 92.11 & 93.50 & 92.05 & 93.46 \\
   & 225 & 92.40 & 93.63 & 92.43 & 93.88 & 92.89 & 94.09 \\
   & 250 & 92.70 & 93.79 & 92.38 & 93.66 & 92.56 & 93.67 \\ \hline
   & 50 & 88.99 & 92.00 & 89.50 & 92.73 & 89.59 & 92.83 \\
   & 75 & 91.10 & 93.57 & 91.14 & 93.59 & 90.91 & 93.76 \\
   & 100 & 90.84 & 93.12 & 92.10 & 94.27 & 91.66 & 93.77 \\
   & 125 & 92.06 & 93.80 & 92.06 & 94.18 & 92.25 & 94.27 \\
  100\% & 150 & 91.95 & 93.89 & 92.56 & 94.65 & 92.60 & 94.59 \\
   & 175 & 92.41 & 94.23 & 92.87 & 94.70 & 92.27 & 94.35 \\
   & 200 & 92.34 & 94.21 & 93.11 & 94.88 & 92.99 & 94.88 \\
   & 225 & 92.73 & 94.31 & 93.25 & 94.85 & 93.06 & 94.89 \\
   & 250 & 93.32 & 94.90 & 92.74 & 94.41 & 92.90 & 94.61 \\
   \hline
\end{tabular}
\caption{Simulated coverage probabilities of Hall-Wellner bands in per cent where $k=10$}
\label{tab:HW10}
\end{table}

\begin{table}[ht]
\centering
\scriptsize
\begin{tabular}{c|c|cc|cc|cc}
  \multicolumn{2}{c|}{set-ups} & \multicolumn{2}{c|}{$N(0,1)$} & \multicolumn{2}{c|}{$Poi(1)-1$}  & \multicolumn{2}{c}{weird}  \\ \hline
  $p$ & $n$ & old & new & old & new & old & new \\ \hline
   & 50 & 88.39 & 89.53 & 89.93 & 91.16 & 89.14 & 90.05 \\
   & 75 & 90.39 & 90.99 & 90.67 & 91.39 & 90.74 & 91.52 \\
   & 100 & 91.39 & 91.89 & 91.75 & 92.26 & 91.92 & 92.30 \\
   & 125 & 91.69 & 92.17 & 92.54 & 93.03 & 92.52 & 92.90 \\
  25\% & 150 & 92.07 & 92.51 & 92.57 & 92.83 & 92.92 & 93.28 \\
   & 175 & 92.24 & 92.55 & 92.76 & 92.97 & 92.8 & 93.07 \\
   & 200 & 93.01 & 93.27 & 92.58 & 93.04 & 92.96 & 93.21 \\
   & 225 & 93.11 & 93.32 & 92.99 & 93.25 & 93.16 & 93.42 \\
   & 250 & 92.88 & 93.03 & 93.60 & 93.69 & 93.05 & 93.31 \\ \hline
   & 50 & 88.77 & 90.02 & 89.28 & 90.66 & 89.53 & 90.89 \\
   & 75 & 90.43 & 91.37 & 91.27 & 92.29 & 90.87 & 91.92 \\
   & 100 & 91.74 & 92.45 & 91.99 & 92.74 & 91.94 & 92.62 \\
   & 125 & 92.12 & 92.82 & 92.49 & 93.11 & 92.08 & 92.61 \\
  50\% & 150 & 92.12 & 92.59 & 93.05 & 93.60 & 92.56 & 93.18 \\
   & 175 & 92.43 & 93.00 & 92.96 & 93.25 & 93.00 & 93.45 \\
   & 200 & 92.94 & 93.22 & 93.08 & 93.50 & 92.69 & 93.28 \\
   & 225 & 92.67 & 93.01 & 93.20 & 93.92 & 93.51 & 93.90 \\
   & 250 & 93.41 & 93.58 & 93.03 & 93.47 & 92.74 & 93.25 \\ \hline
   & 50 & 89.17 & 90.83 & 89.70 & 91.53 & 89.65 & 91.36 \\
   & 75 & 90.99 & 92.23 & 91.29 & 92.41 & 91.21 & 92.32 \\
   & 100 & 91.15 & 92.12 & 92.36 & 93.24 & 92.26 & 93.37 \\
   & 125 & 92.21 & 93.15 & 92.57 & 93.44 & 92.49 & 93.50 \\
  75\% & 150 & 92.64 & 93.39 & 93.08 & 93.98 & 92.68 & 93.46 \\
   & 175 & 92.50 & 93.21 & 93.25 & 93.91 & 92.67 & 93.44 \\
   & 200 & 92.73 & 93.39 & 93.34 & 94.15 & 93.10 & 93.78 \\
   & 225 & 92.91 & 93.51 & 93.18 & 93.96 & 93.61 & 94.11 \\
   & 250 & 93.34 & 94.01 & 93.84 & 94.39 & 92.92 & 93.67 \\ \hline
   & 50 & 90.11 & 92.25 & 90.74 & 92.78 & 90.52 & 92.90 \\
   & 75 & 91.22 & 92.70 & 92.57 & 94.08 & 91.94 & 93.44 \\
   & 100 & 92.01 & 93.62 & 93.09 & 94.37 & 92.49 & 93.84 \\
   & 125 & 92.39 & 93.45 & 93.28 & 94.64 & 93.59 & 94.56 \\
  100\% & 150 & 92.77 & 94.10 & 93.47 & 94.59 & 93.28 & 94.41 \\
   & 175 & 93.22 & 94.35 & 93.66 & 94.84 & 94.00 & 94.93 \\
   & 200 & 93.88 & 94.73 & 93.79 & 94.78 & 93.90 & 94.84 \\
   & 225 & 93.74 & 94.62 & 94.25 & 95.07 & 93.69 & 94.70 \\
   & 250 & 93.61 & 94.56 & 93.94 & 94.95 & 94.17 & 95.12 \\
   \hline
\end{tabular}
\caption{Simulated coverage probabilities of Hall-Wellner bands in per cent where $k=20$}
\label{tab:HW20}
\end{table}

\clearpage

\begin{table}[ht]
\centering
\scriptsize
\begin{tabular}{c|cc|cc|cc|cc|cc|cc}
  & \multicolumn{6}{|c|}{Hall-Wellner} &  \multicolumn{6}{|c}{Equal Precision} \\ \hline
  & \multicolumn{2}{c|}{$N(0,1)$} & \multicolumn{2}{c|}{$Poi(1)-1$}  & \multicolumn{2}{c|}{weird} & \multicolumn{2}{c|}{$N(0,1)$} & \multicolumn{2}{c|}{$Poi(1)-1$}  & \multicolumn{2}{c}{weird} \\ \hline
  $n$ & old & new & old & new & old & new & old & new & old & new & old & new \\ \hline
  50 & 87.61 & 88.64 & 88.80 & 90.03 & 88.72 & 89.75 & 87.08 & 88.06 & 89.55 & 90.36 & 89.43 & 90.48 \\
  75 & 90.15 & 90.75 & 90.71 & 91.46 & 90.74 & 91.48 & 89.36 & 90.21 & 91.79 & 92.33 & 91.62 & 92.14 \\
  100 & 90.78 & 91.31 & 91.13 & 91.58 & 91.77 & 92.23 & 90.74 & 91.24 & 93.20 & 93.43 & 92.43 & 92.72 \\
  125 & 91.72 & 92.07 & 92.02 & 92.31 & 91.58 & 92.23 & 91.35 & 91.77 & 92.96 & 93.37 & 93.19 & 93.54 \\
  150 & 92.29 & 92.54 & 92.60 & 92.85 & 92.90 & 93.21 & 91.33 & 91.52 & 93.64 & 93.85 & 93.51 & 93.76 \\
  175 & 92.27 & 92.53 & 92.48 & 92.83 & 92.86 & 93.17 & 92.45 & 92.64 & 93.91 & 94.14 & 94.06 & 94.39 \\
  200 & 92.67 & 92.85 & 92.89 & 93.02 & 93.18 & 93.34 & 92.87 & 93.04 & 93.68 & 93.80 & 93.47 & 93.83 \\
  225 & 92.98 & 93.05 & 92.60 & 92.71 & 93.26 & 93.68 & 92.69 & 92.93 & 93.85 & 94.12 & 94.16 & 94.22 \\
  250 & 92.78 & 92.99 & 92.93 & 93.14 & 93.14 & 93.32 & 93.20 & 93.43 & 93.68 & 93.81 & 94.25 & 94.30 \\
   \hline
\end{tabular}
\caption{Simulated coverage probabilities of confidence bands for $F_1$ in per cent where $p=0$}
\label{tab:HWEPcont}
\end{table}


\section{Real Data Example}
\label{sec:data_ex}

We applied the present discretization adjustment to the \verb=sir.adm= data-set of the \verb=R= package {\it mvna}.
It consists of competing risks data of patients who are in an intensive care unit (ICU),
where the event of primary interest, ``alive discharge out of ICU'', competes against the secondary event ``death in ICU''.
For seeing the difference between the common and the new approach more clearly,
we analyzed the subset of all male patients suffering from pneumonia.
Out of these $n=63$ individuals, five have been right-censored and 41 out of all 44 type 1 events fell into the time interval $[5,55]$, which we chose for confidence band construction.
Due to the worse performance of the wild bootstrap based on standard normal multipliers as seen in Section~\ref{sec:simus},
we derived these bands by only using centered unit Poisson variates.
In order to minimize the computational error in the quantile-finding process,
99,999 wild bootstrap iterations have been conducted.
The confidence bands resulting from the weird bootstrap almost coincide with those just described.
Therefore, they are not shown.

\begin{figure}[ht]
 \includegraphics[width=1.0\textwidth, height=0.4\textwidth]{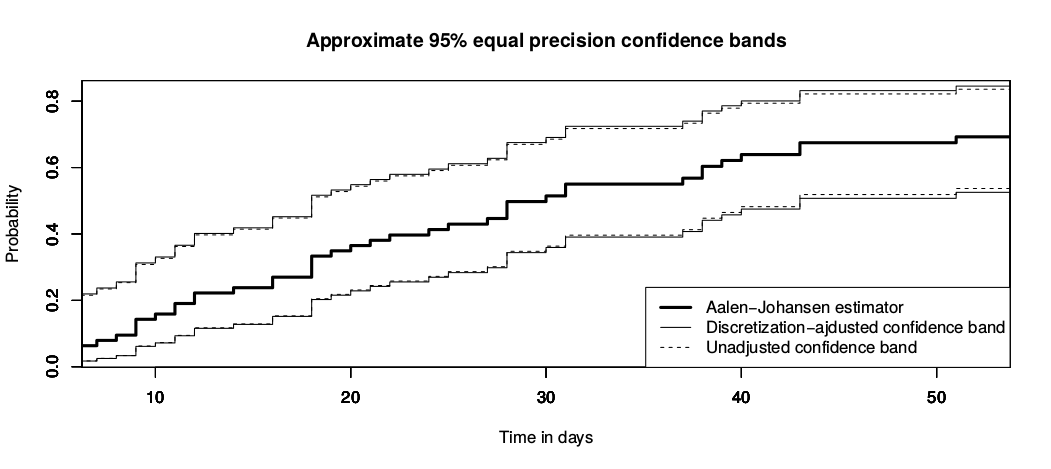} \\
 \includegraphics[width=1.0\textwidth, height=0.4\textwidth]{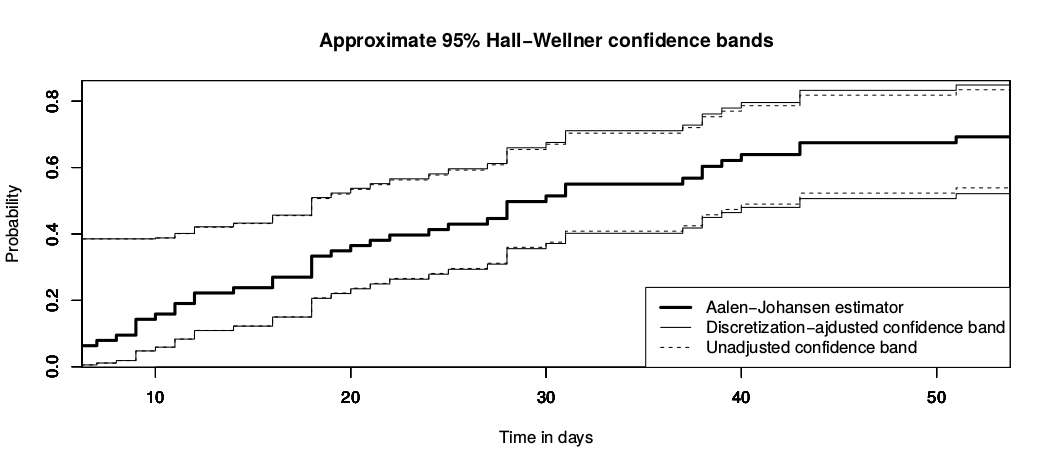}
  \caption{Asymptotic $95\%$ equal precision (upper) and Hall-Wellner bands (lower panel) for the cumulative incidence function of the competing risk ``alive discharge out of ICU'' for male patients suffering from pneumonia.}
 \label{fig:bands}
\end{figure}

The resulting confidence bands are presented in Figure~\ref{fig:bands}. 
For both kinds of bands,  (equal precision bands in the upper panel, Hall-Wellner bands in the lower panel),
we see that the discretization adjustment leads to a widening in comparison to the unadjusted bands.
This is in line with the results of the simulation study in Section~\ref{sec:simus},
where the unadjusted bands appeared to be the most liberal, i.e. the smallest.
In particular, the adjusted equal precision bands cover an additional area of $2.1$ percentage points at the terminal time point $t=55$,
whereas this deviance even amounts to $3.3$ percentage points for the Hall-Wellner bands.
This might not appear to be much at a first glance at the plots in Figure~\ref{fig:bands}. 
But in fact, it may be the cause for a formidable improvement of the bands' coverage probability:
The simulation results of Section~\ref{sec:simus} for $k=20$, discretization probability $p=100\%$, and sample sizes $n \in \{50, 75\}$ suggest that the adjusted wild bootstrap procedure
might improve the coverage probabilities of both kinds of bands by approximately two percentage points.
With a view towards the liberal behaviour of the unadjusted bands,
these enhancements of the coverage probabilities are highly worthwhile.

\section{Discussion and Future Research}
\label{sec:disc}

In this article, we analyzed a discontinuity adjustment of the common wild bootstrap applied to right-censored competing risks data.
This adjustment is absolutely necessary, as ties in the data introduce an asymptotic dependence between multiple cause-specific Nelson-Aalen estimators.
The common wild bootstrap fails in reproducing this effect since it establishes independence for all sample sizes.
The problem is even more involved for Aalen-Johansen estimators of cumulative incidence functions,
which are non-linear functionals of all cause-specific hazards.
Simulation results reported the striking liberality of the unadjusted bands which also fail to keep the nominal level asymptotically.
Instead, the discretization-adjusted wild bootstrap greatly improves the coverage probability.
This effect is more pronounced the more discrete the event times are.
But even in the absolutely continuous case, the suggested procedure appears to perform preferably.
Therefore, we advise to always use the adjustment when right-censored competing risks data shall be analyzed.
The real data example reveals that the discontinuity adjustment does only lead to slight widening of the common wild bootstrap-based bands
which is already enough to improve the coverage accuracy greatly.

The presented wild bootstrap approach may be extended to more general Markovian multi-state models since the martingale arguments of Appendix~\ref{sec:app_wbs_univ_nae} still apply.
A still open question is, whether the common wild bootstrap also fails in case of tied survival data which are assumed to follow the Cox proportional hazards model \citep{cox72}.
See \cite{lin93} for the wild bootstrap applied to absolutely continuous survival data following the Cox model.
And if it fails, does the method proposed in this article require further modification?

\section*{Acknowledgements}

The author likes to thank Markus Pauly for helpful discussions.
He also appreciates the support by the German Research Foundation (DFG).

\appendix

\section{Detailed Derivation of the Asymptotic Variance of the Nelson-Aalen Estimator}
\label{sec:app_asy_var}

Define $H^{uc}(t) = \sum_{j=1}^k E N_{j1}(t)$ as the probability that an uncensored event due to any cause occurs until time $t$.
According to \cite{vaart96}, p. 383f, we have $\sqrt{n} (\wh A - A ) \oDo \int_0^\cdot M^{uc}(\d u) /\bar H(u)$,
where $M^{uc}(t) = G^{uc}(t) - \int_0^t \bar {G}(u) \d A(u)$ is a zero-mean Gaussian martingale.
Its variance function is determined by
$$ E G^{uc}(s) G^{uc}(t) = H^{uc}(s\wedge t) - H^{uc}(s) H^{uc}(t),$$
$$ E \bar{G}(s) \bar {G}(t) = \bar H(s \vee t) - \bar H(s) \bar H(t),$$
$$ E G^{uc}(s) \bar{G}(t) = (H^{uc}(s) - H^{uc}(t-)) 1\{t \leq s\} - H^{uc}(s) \bar H (t).$$
Note that $A(t) = \int_0^t H^{uc}(\d u) / \bar H(u)$.
Thus, for $s \leq t$, the covariance function of $M^{uc}$ at $(s,t)$ is
\begin{align*}
 & E(M^{uc}(s) M^{uc}(t)) \\
 & = H^{uc}(s) - H^{uc}(s) H^{uc}(t)
  + \int_0^s \int_0^t \frac{\bar H(u \vee v) - \bar H(u) \bar H(v)}{\bar H(u) \bar H(v)} \d H^{uc}(v) \d H^{uc}(u) \\
 & - \int_0^t [(H^{uc}(s) - H^{uc}(u-)) 1\{u \leq s\} - H^{uc}(s) \bar H (u)] \frac{\d H^{uc}(u)}{\bar H(u)} \\
 & - \int_0^s [(H^{uc}(t) - H^{uc}(v-)) 1\{v \leq t\} - H^{uc}(t) \bar H (v)] \frac{\d H^{uc}(v)}{\bar H(v)} \\
 & = H^{uc}(s) - H^{uc}(s) H^{uc}(t)
 + \int_0^s \int_u^t \Big[\frac{1}{\bar H(u)} - 1 \Big] \d H^{uc}(v) \d H^{uc}(u)
 + \int_0^s \int_0^u \Big[\frac{1}{\bar H(v)} - 1 \Big] \d H^{uc}(v) \d H^{uc}(u) \\
 & - (H^{uc}(s) + H^{uc}(t)) A(s) + 2 \int_0^s H^{uc}(u-) \d A(u) + 2 H^{uc}(s) H^{uc}(t) \\
 & = H^{uc}(s) + H^{uc}(s) H^{uc}(t) + \int_0^s (H^{uc}(t) - H^{uc}(u)) \Big[\frac{1}{\bar H(u)} - 1 \Big]  \d H^{uc}(u)
 + \int_0^s (A(u) - H^{uc}(u)) \d H^{uc}(u) \\
 & - (H^{uc}(s) + H^{uc}(t)) A(s) + 2 \int_0^s H^{uc}(u-) \d A(u) \\
 & = H^{uc}(s) + H^{uc}(s) H^{uc}(t) + H^{uc}(t) A(s) - H^{uc}(t) H^{uc}(s) - \int_0^s H^{uc}(u) \d A(u)+ \int_0^s H^{uc}(u) \d H^{uc}(u) \\
 & + \int_0^s (A(u) - H^{uc}(u)) \d H^{uc}(u) - (H^{uc}(s) + H^{uc}(t)) A(s) + 2 \int_0^s H^{uc}(u-) \d A(u) \\
 & = H^{uc}(s) - H^{uc}(s) A(s) - \int_0^s H^{uc}(u) \d A(u) + \int_0^s A(u) \d H^{uc}(u) + 2 \int_0^s H^{uc}(u-) \d A(u) \\
 & = H^{uc}(s) - \int_0^s H^{uc}(u) \d A(u) + \int_0^s H^{uc}(u-) \d A(u) \\
 & = \int_0^s \bar H(u) \d A(u) - \int_0^s \Delta H^{uc}(u) \d A(u)
 = \int_0^s \bar H(u) (1 - \Delta A(u)) \d A(u)
\end{align*}
We conclude, as in \cite{vaart96}, that
$$\sqrt{n} (\wh A - A ) \oDo \int_0^\cdot \frac{1}{\bar H(u)} \d M^{uc}(u) \sim Gauss \Big(0, \int_0^\cdot \frac{1 - \Delta A(u)}{\bar H(u)} \d A(u) \Big),$$
where $Gauss$ again indicates that the limit process is a Gaussian martingale.

The very same calculations hold true if each $H^{uc}$ is replaced with $H^{uc}_1$ for the subdistribution function of an uncensored type 1 event.
Therefore, we have for Nelson-Aalen estimators for cause-specific cumulative hazard functions that
$$\sqrt{n} (\wh A_1 - A_1 ) \oDo \int_0^\cdot \frac{1}{\bar H(u)} \d \mb M_1^{uc}(u) \sim Gauss \Big(0, \sigma^2_1(\cdot) \Big).$$

\section{Consistency of the Wild Bootstrap for the Univariate Nelson-Aalen Estimator}

\label{sec:app_wbs_univ_nae}

\begin{proof}
 Without loss of generality, assume that $0, K \in D_1$ for simplifying notation.
 Write $D_1 = \{ d_0, d_1, \dots, d_J \}$ with the natural ordering $d_j < d_{j+1}$ for all $j=1, \dots, J$.
 Then $[0,K] \setminus D_1 = \bigcup\limits_{j=1}^J (d_{j-1}, d_j)$.
 As argued in \cite{bluhmki15}, it is now straightforward to show
 that each process $(\wh W_1(t) - \wh W_1(d_{j-1}))_t$ on each interval $[d_{j-1}, d_j)$, $j=1, \dots, J$,
 defines a square-integrable martingale.
 Since such martingales can be extended to the right boundary of the time interval,
 we may \emph{define} the boundary values
 $\wh W_1(d_j) - \wh W_1(d_{j-1}) := \lim_{t \uparrow d_j} \wh W_1(t) - \wh W_1(d_{j-1})$,
 and this procedure introduces square-integrable martingales on the whole $[d_{j-1}, d_j]$.

 First, we notice that the conditional weak convergence in probability of the processes $(\wh W_1(t) - \wh W_1(d_{j-1}))_t$ on each interval $[d_{j-1}, d_j]$, $j=1, \dots, J$,
 is already implied by exactly the same Rebolledo's martingale central limit theorem arguments as in \cite{bluhmki15}.
 Denote the limit Gaussian martingale processes as $(\tilde U_{1j}(t))_{t \in [d_{j-1}, d_j]}$, $j=1,\dots, J$.
 Due to the martingale extension above,
 Rebolledo's limit theorem implies the almost sure continuity of $\tilde U_{1j}$ on each time interval.
 Furthermore, these are zero-mean processes with variance function $t \mapsto \sigma^2_1(t) - \sigma^2_1(d_{j-1})$.

 Due to the continuity of the limit processes $\tilde U_{1j}$ on the intervals $[d_{j-1}, d_j]$,
 we are able to switch from the Skorohod topology to the more convenient $\sup$-norm metrization;
 see the discussion in Section~II.8 in \cite{abgk93}.
 At each $t = d_j$, the weak conditional convergence in distribution of $\Delta \wh W_1(d_j)$ holds in probability by the already argued convergence of all finite-dimensional conditional distributions.
 Therefore, the independence of the (bootstrapped) Nelson-Aalen increments imply
 that, as $\nri$, the conditional distribution of
 $$ \Big(\Delta \wh W_1(d_0), \wh W_1(t_0) - \wh W_1(d_0), \Delta \wh W_1(d_1), \dots, \wh W_1(t_{J-1}) - \wh W_1(d_{J-1}), \Delta \wh W_1(d_J) \Big)_{t_0 \in [d_0, d_1], \dots, t_{J-1} \in [d_{J-1}, d_J]} $$
 given $\mac F_0$ converges weakly in probability to the distribution of
 $$ (\tilde U_1(d_0), \tilde U_{11}(t_1), \tilde U_1(d_1), \dots, \tilde U_{1J}(t_{J}), \tilde U_1(d_J))_{t_1 \in [d_0, d_1], \dots, t_{J} \in [d_{J-1}, d_J]} $$
 on the product Space $\R \times D[d_0,d_1] \times \R \times \dots \times D[d_{J-1}, d_J] \times \R$ equipped with the $\sup$-$\max$-norm.
 Here, all components are independent,
 and the normally distributed random variables $\tilde U_1(d_j)$ have mean zero and variance
 $\Delta \sigma^2_1(d_j)$, $j=1,\dots,J$.

 Applying the functional, which is continuous with respect to the $\max$-$\sup$-norm,
 $$\psi: \R \times D[d_0,d_1] \times \R \times \dots \times D[d_{J-1}, d_J] \times \R
  \longrightarrow  D[d_0,d_1] \times \dots \times D[d_{J-1}, d_J] \times \R,$$
  \begin{align*}
  & (x_0, y_1(t_1), x_1, \dots, y_J(t_J), x_J)_{t_1, \dots, t_J}
  \longmapsto \Big(x_0 + y_1(t_1), \ x_0 + y_1(d_1) + x_1 + y_2(t_2), \\
  & \dots,
  x_0 + \sum_{j=1}^{J-1} (y_j(d_j) + x_j) + y_{J}(t_J), \ x_0 + \sum_{j=1}^{J} (y_j(d_j) + x_j) \Big)_{t_1, \dots, t_J}
  \end{align*}
 to the previous limit theorem, the continuous mapping theorem implies that:
 Given $\mac F_0$ and as $\nri$,
 the conditional distribution of
 $$ \Big(\wh W_1(t_1), \dots, \wh W_1(t_{J}), \wh W_1(d_J) \Big)_{t_1 \in [d_0, d_1], \dots, t_{J} \in [d_{J-1}, d_J]} $$
 converges weakly (on the product-function space $D[d_0,d_1] \times \dots \times D[d_{J-1}, d_J] \times \R$ equipped with the $\max$-$\sup$-norm) in probability to the distribution of
 $$ \Big(U_1(t_1), \dots, U_1(t_{J}), U_1(d_J) \Big)_{t_1 \in [d_0, d_1], \dots, t_{J} \in [d_{J-1}, d_J]}. $$
 Here the right boundary values are again considered as the left-hand limits given by the martingale extension theorem.
 The process $(U_1(t))_{t \in [0,K]}$ is a zero-mean Gaussian martingale with variance function
 $t \mapsto \sigma^2_1(t)$ with, in general, discontinuous sample paths.

 Finally, we apply the continuous functional
 \begin{align*}
  \phi : \ & D[d_0,d_1] \times \dots \times D[d_{J-1}, d_J] \times \R \longrightarrow D[0,K], \\
  & (y_1(t_1), \dots, y_J(t_J), x_J)_{t_1, \dots, t_J} \longmapsto \Big( \sum_{j=1}^J y_j(t) \cdot 1_{[d_{j-1}, d_j)} (t) + x_J \cdot 1_{\{K\}} (t) \Big)_{t \in [0,K]}
 \end{align*}
 in order to obtain the desired conditional weak convergence for $(\wh W_1(t))_{t \in [0,K]}$ in probability.
\end{proof}

\section{Derivation of the Asymptotic Covariance of Multiple Nelson-Aalen Estimators for Cumulative Cause-Specific Hazards}

\label{sec:app_asy_cov_csh}

Consider the situation of two competing risks.
In order to derive the asymptotic covariance of two cause-specific Nelson-Aalen estimators, we first note that, as $\nri$,
$$\sqrt{n} (\wh A - A) = \sum_{j = 1}^2 \sqrt{n}(\wh A_j - A_j) \oDo U_1 + U_2 = U$$
by the continuous mapping theorem.
The covariance function of $U$ is given by
\begin{align*}
 cov(U(s), U(t)) = \sigma^2(s \wedge t) = \int_0^{s \wedge t} \frac{1 - \Delta A(u)}{\bar H(u)} \d A(u)
\end{align*}
but on the other hand
\begin{align*}
 cov(U(s), U(t)) = \sum_{j=1}^2 cov(U_j(s), U_j(t)) + cov(U_1(s),U_2(t)) + cov(U_1(t), U_2(s)).
\end{align*}
Solving for the unknown covariances on the right-hand side of the previous display, we obtain
\begin{align*}
 cov(U_1(s),U_2(t)) & + cov(U_1(t), U_2(s))
  = \int_0^{s \wedge t} \frac{1 - \Delta A(u)}{\bar H(u)} \d A(u)
 - \sum_{j=1}^2 \int_0^{s \wedge t} \frac{1 - \Delta A_j(u)}{\bar H(u)} \d A_j(u) \\
 & = - \int_0^{s \wedge t} \frac{\Delta A_1(u)}{\bar H(u)} \d A_2(u)
 -  \int_0^{s \wedge t} \frac{\Delta A_2(u)}{\bar H(u)} \d A_1(u)
\end{align*}
Due to symmetry and inductively, it follows that
$$ \sigma_{j \ell}(s \wedge t) =  cov(U_j(s),U_\ell(t)) = - \int_0^{s \wedge t} \frac{\Delta A_j(u)}{\bar H(u)} \d A_\ell(u) $$
for $j \neq \ell, \ j,\ell = 1, \dots, k,$ even in the situation of $k \in \N$ competing risks.

\section{Consistency of the Discretization-Adjusted Wild Bootstrap for the Multivariate Nelson-Aalen Estimator}

\label{sec:app_wbs_mult_nae}

The proof of tightness follows along the same lines as that of Theorem~\ref{thm:main}.
It only remains to calculate the finite-dimensional marginal limit distributions.
These are calculated with the help of Theorem~A.1 in \cite{beyersmann13}:
Therefore, we abbreviate
$$\wh W_j(t) = \sum_{i=1}^n \xi_{jji} Z_{n,jji}
  + \sum_{\ell=1}^k sign(\ell - j) \sum_{i=1}^n \xi_{j \ell i} \tilde Z_{n,j \ell i}
  + \sum_{\ell=1}^k sign(\ell - j) \sum_{i=1}^n \xi_{\ell ji} \tilde Z_{n,\ell ji}. $$
Further, consider arbitrary points of time $0 \leq t_1 \leq \dots \leq t_m \leq K$
and the vector
$$(\wh W_1(t_1), \dots, \wh W_1(t_m), \wh W_2(t_1), \dots, \wh W_2(t_m) \dots, \wh W_k(t_1), \dots, \wh W_k(t_m)). $$
Due to analogy, it is enough to calculate the entries corresponding to $(\wh W_1(t_1), \wh W_1(t_2))$ and $(\wh W_1(t_1), \wh W_2(t_2))$ of the limit of the matrix
\begin{align*}
 \wh \Gamma := & \Big( 1\{j = \tilde j \}  \sum_{i=1}^n \Big[  Z_{n,jji}(t_a) Z_{n, jj i}(t_b)
  + \sum_{\ell \neq j} \tilde Z_{n,j \ell i}(t_a) \tilde Z_{n, j \ell i}(t_b)
  + \sum_{\ell \neq j} \tilde Z_{n,\ell j i}(t_a) \tilde Z_{n, \ell j i}(t_b) \Big] \\
   & - 1\{j \neq \tilde j \} \sum_{i=1}^n \Big[ \tilde Z_{n,j \tilde j i}(t_a) \tilde Z_{n, j \tilde j i}(t_b)
  + \tilde Z_{n, \tilde j j i}(t_a) \tilde Z_{n, \tilde j j i}(t_b) \Big]
  \Big)_{a,b=1, \dots, m; \ j, \tilde j = 1 , \dots, k}.
\end{align*}
We start by calculating the entry for $a=1, b=2, j=\tilde j=1$, that is,
\begin{align*}
  & \sum_{i=1}^n \Big[ Z_{n,11i}(t_1) Z_{n, 11i}(t_2)
    + \sum_{\ell = 2}^k \tilde Z_{n,1\ell i}(t_1) \tilde Z_{n, 1\ell i}(t_2)
    + \sum_{\ell = 2}^k \tilde Z_{n,\ell 1 i}(t_1) \tilde Z_{n, \ell 1 i}(t_2) \Big] \\
  & = n \sum_{i=1}^n  \int_0^{t_1} {\frac{Y(u) - \Delta N(u)}{Y(u)}} \frac{\d N_{1i}(u)}{Y^2(u)}
  + \frac12 \sum_{\ell = 2}^k \Big[ n \sum_{i=1}^n \int_0^{t_1} {\frac{\Delta N_1(u)}{Y(u)}} \frac{\d N_{\ell i}(u)}{Y^2(u)}
  +  n \sum_{i=1}^n \int_0^{t_1} {\frac{\Delta N_\ell (u)}{Y(u)}} \frac{\d N_{1 i}(u)}{Y^2(u)} \Big] \\
  & = n \int_0^{t_1} {\frac{Y(u) - \Delta N(u)}{Y^2(u)}} \frac{\d N_{1}(u)}{Y(u)}
  + \frac12 \sum_{\ell = 2}^k \Big[ n \int_0^{t_1} {\frac{\Delta N_1(u)}{Y^2(u)}} \frac{\d N_{\ell}(u)}{Y(u)}
  +  n \int_0^{t_1} {\frac{\Delta N_\ell (u)}{Y^2(u)}} \frac{\d N_{1}(u)}{Y(u)} \Big] \\
  & = n \int_0^{t_1} {\frac{Y(u) - \Delta N_1(u)}{Y^2(u)}} \frac{\d N_{1}(u)}{Y(u)}.
\end{align*}
Here, the last equality follows from $\Delta N_1 \d N_{\ell} = \Delta N_\ell \d N_{1} $ for all $\ell$.
By the Glivenko-Cantelli theorem in combination with the continuous mapping theorem,
it follows that the quantity in the previous display converges to $\sigma_1^2(t_1)$ in probability as $\nri$.

Now, consider the entry of $\wh \Gamma$ for $a=1, b=2, j=\tilde j=2$:
\begin{align*}
 - \sum_{i=1}^n \Big[ \tilde Z_{n,12i}(t_1) \tilde Z_{n,12i}(t_2) + \tilde Z_{n,21i}(t_1) \tilde Z_{n,21i}(t_2) \Big]
 = - \frac12 \Big[ n \int_0^{t_1} {\frac{\Delta N_1 (u)}{Y^2(u)}} \frac{\d N_{2}(u)}{Y(u)}
  + n \int_0^{t_1} {\frac{\Delta N_2 (u)}{Y^2(u)}} \frac{\d N_{1}(u)}{Y(u)} \Big].
\end{align*}
By the same arguments as before, this is a consistent estimator for $\sigma_{12}(t_1)$ as $\nri$.

\section{Derivation of the Asymptotic Covariance Function of an Aalen-Johansen Estimator in the Competing Risks Set-Up}.

\label{sec:app_asy_var_aje}

In order to derive the the covariance function of
\begin{align*}
 \int_0^\cdot \frac{1 - F_2(u-) - F_1( \cdot )}{1 - \Delta A(u)} \d U_1(u)
 + \int_0^\cdot \frac{F_1(u-) - F_1(\cdot)}{1 - \Delta A(u)} \d U_2(u),
\end{align*}
at any $(s,t) \in [0,K]^2$, we exemplarily calculate the covariance function of the first integral and the covariance function between both integrals.
Hence, as the covariance function $(s,t) \mapsto \sigma^2_1(s \wedge t)$ of $U_1$ only increases along the diagonal,
\begin{align*}
 & cov \Big(\int_0^s \frac{1 - F_2(u-) - F_1( s )}{1 - \Delta A(u)} \d U_1(u) , \int_0^t \frac{1 - F_2(u-) - F_1( t )}{1 - \Delta A(u)} \d U_1(u) \Big) \\
 & = \int_0^{s \wedge t} \frac{(1 - F_2(u-) - F_1( s )) (1 - F_2(u-) - F_1( t ))}{(1 - \Delta A(u))^2} \d \sigma^2_1(u) \\
 & = \int_0^{s \wedge t} \frac{(1 - F_2(u-) - F_1(s))(1 - F_2(u-) - F_1(t))}{\bar H(u)} \frac{1- \Delta A_1(u)}{(1- \Delta A(u))^2} \d A_1(u).
\end{align*}
Furthermore, we similarly have for the covariance between both integrals that
\begin{align*}
 & cov \Big(\int_0^s \frac{1 - F_2(u-) - F_1( s )}{1 - \Delta A(u)} \d U_1(u) , \int_0^t \frac{F_1(u-) - F_1( t )}{1 - \Delta A(u)} \d U_2(u) \Big) \\
 & = \int_0^{s \wedge t} \frac{(1 - F_2(u-) - F_1( s )) (F_1(u-) - F_1( t ))}{(1 - \Delta A(u))^2} \d \sigma_{12}(u) \\
 & = - \int_0^{s \wedge t} \frac{(1 - F_2(u-) - F_1(s))(F_1(u-) - F_1(t))}{\bar H(u)} \frac{\Delta A_1(u)}{(1- \Delta A(u))^2} \d A_2(u).
\end{align*}
Finally, including also the remaining two analogous terms, we obtain the following asymptotic covariance function of the Aalen-Johansen estimator
for the first cumulative incidence function
as the sum of all four covariance functions:
\begin{align*}
 (s,t) \mapsto & \int_0^{s \wedge t} \frac{(1 - F_2(u-) - F_1(s))(1 - F_2(u-) - F_1(t))}{\bar H(u)} \frac{1- \Delta A_1(u)}{(1- \Delta A(u))^2} \d A_1(u) \\
    & + \int_0^{s \wedge t} \frac{(F_1(u-) - F_1(s))(F_1(u-) - F_1(t))}{\bar H(u)} \frac{1- \Delta A_2(u)}{(1- \Delta A(u))^2} \d A_2(u) \\
    & - \int_0^{s \wedge t} \frac{(1 - F_2(u-) - F_1(s))(F_1(u-) - F_1(t))}{\bar H(u)} \frac{\Delta A_1(u)}{(1- \Delta A(u))^2} \d A_2(u) \\
    & - \int_0^{s \wedge t} \frac{(1 - F_2(u-) - F_1(t))(F_1(u-) - F_1(s))}{\bar H(u)} \frac{\Delta A_2(u)}{(1- \Delta A(u))^2} \d A_1(u) .
\end{align*}


\renewcommand{\refname}{References}
\bibliography{literatur}
\bibliographystyle{plainnat}

\newpage

\section*{Erratum}





{\centering
Dennis Dobler \\
{\small TU Dortmund University, \ \ Department of Statistics, \\
and University Alliance Ruhr, \ \ Research Center Trustworthy Data Science and Security,  \\
Joseph-von-Fraunhofer-Str.\ 25, \ \  44227 Dortmund, \ \ Germany \\
email: dobler@statistik.tu-dortmund.de \ \ \\[0.5cm]
and\\[0.5cm]
Merle Munko \\
{\small Otto-von-Guericke University Magdeburg, \ \ Department of Mathematics, \\
Postfach 4120, \ \ 39106 Magdeburg, \ \ Germany \\
email: merle.munko@ovgu.de \ \ \\
}}

\vspace{0.5cm}

{
\centering

\date{\today}

}

\vspace{0.5cm}

\begin{center}
{\bf Summary} \vspace{-.4cm}
\end{center}
\justifying
\noindent
The wild bootstrap is the resampling method of choice in survival analytic applications.
Theoretic justifications rely on the assumption of existing intensity functions
which is equivalent to an exclusion of ties among the event times.
However, such ties are omnipresent in practical studies.
It turns out that the wild bootstrap should only be applied in a modified manner
that corrects for altered limit variances and emerging dependencies.
This again ensures the asymptotic exactness of inferential procedures.
An analogous necessity is the use of the Greenwood-type variance estimator for Nelson-Aalen estimators
which is particularly preferred in tied data regimes.
All theoretic arguments are transferred to bootstrapping Aalen-Johansen estimators for cumulative incidence functions in competing risks.
An extensive simulation study as well as an application to real competing risks data of male intensive care unit patients suffering from pneumonia illustrate the practicability of the proposed technique.

\noindent{\bf Keywords:}
Aalen-Johansen estimator;
Counting process;
Discontinuous cumulative hazard functions;
Discontinuous cumulative incidence functions;
Greenwood-type variance estimator;
Nelson-Aalen estimator;
Survival analysis;
Tied event times;
Wild bootstrap.

\newpage

We consider the same competing risks setup as above, i.e., we assume that there are $k\in\mathbb N$ competing risks and $n\in\mathbb N$ i.i.d.\ random event times $T_1,...,T_n$, which are independently right-censored and distributed as a random variable $T \sim S$. Here, $S$ denotes the survival function, i.e., $S(t) = P(T>t)$ for all $t \geq 0$; $S$ need not be continuous. Then, we denote the probability that an individual is under observation at time $t-$, that is, \emph{just before} time $t$, by $\bar{H}(t) = P(\min(T,C) \geq t) = S(t-)G(t-)$ for all $t \geq 0$.
Here, $C\sim G$ with survival function $G(t) = P(C>t)$ denotes a generic censoring time which is assumed to be independent of $T$.
Also, $t \mapsto {f}(t-)$ denotes the left-continuous version of a right-continuous function $f: [0,\infty) \to \mathbb{R}$.
Furthermore, let $\widehat{A}_j$ denote the cause-specific Nelson-Aalen estimator for the cumulative hazard function $A_j$ of type $j$ events, $\widehat{S}$ the Kaplan-Meier estimator for the Survival function $S$, and $\widehat{F}_j = \int_0^. \widehat{S}(u-)\mathrm{d}\widehat{A}_j(u)$ the Aalen-Johansen estimator for the cumulative incidence function ${F}_j = \int_0^. {S}(u-)\mathrm{d}{A}_j(u)$ for all $j\in\{1,...,k\}$. In addition to the assumptions in the main manuscript, it is actually required that $\bar{H}(K)>0$ for $K \geq 0$ to ensure finite variances $\sigma_j^2(K), j\in\{1,\dots,k\},$ in Theorem~3.1 above.


Theorem~4.1 above states for $k=2$ competing risks that
\begin{align*}
    \sqrt{n}(\widehat{F}_1 - F_1) \xrightarrow{d} U_{F_1}
\end{align*} as $n\to\infty$ on the càdlàg space $D[0,K]$ equipped with the sup-norm, where $U_{F_1}$ is a zero-mean Gaussian process with covariance function
\begin{align*}
    \sigma^2_{F_1}: (s,t) \mapsto &\int_0^{s \wedge t} \frac{(1 - F_2(u-) - F_1(s))(1 - F_2(u-) - F_1(t))}{\bar{H}(u)}\frac{\mathrm{d}A_1(u)}{1 - \Delta A(u)}\\
    &+ \int_0^{s \wedge t} \frac{(F_1(u-) - F_1(s))(F_1(u-) - F_1(t))}{\bar{H}(u)}\frac{\mathrm{d}A_2(u)}{1 - \Delta A(u)}\\
    &+ \sum_{u \in D, u\leq s,t} \frac{S^2(u-)}{\bar{H}(u)}\frac{\Delta A_1(u)\Delta A_2(u)}{(1 - \Delta A(u))^2},
\end{align*} where $A = \sum_{j=1}^k A_j$ and $D = \{t\in [0,K] : \Delta A(t)  > 0 \}$ is the set of discontinuity time points.
However, we found that the right-continuous versions $F_1,F_2,S$ must appear in the covariance function above in all occurrences of 
$F_1(u-),F_2(u-),S(u-)$. 

In order to prove this, we go one step back:
By Theorem~3.1 above,  \begin{align*}
    \sqrt{n} \left( \widehat{A}_{1} - {A}_{1}, ...,  \widehat{A}_{k} - {A}_{k} \right) \xrightarrow{d} \left(U_{1}, ..., U_{k}  \right)
\end{align*} holds as $n\to\infty$ on the product space $D^k[0,K]$ equipped with the $\max$-$\sup$ norm,
where $U_{1}, ..., U_{k}$ are zero-mean Gaussian-martingales with 
\begin{align*}
    &\Cov (U_{j}(t), U_{j}(s)) = \int_0^{t \wedge s} \frac{1 - \Delta {A}_{j}(u)}{\bar{H}(u)} \;\mathrm{d}{A}_{j}(u) =: {\sigma}_{j}^2(t \wedge s)
   ,\\& \Cov (U_{j}(t), U_{\ell}(s)) = - \int_0^{t \wedge s} \frac{\Delta {A}_{\ell}(u)}{\bar{H}(u)} \;\mathrm{d}{A}_{j}(u) =: {\sigma}_{j\ell}(t\wedge s)
\end{align*}
 for all $t,s\in [0,K], j,\ell\in\{1,...,k\},j \neq \ell$. We further note that the limit $\left(U_{1}, ..., U_{k}  \right)$ is separable since $G_1^{uc},...,G_k^{uc}$ and $\overline{G}$ in Appendix~A are tight, which follows by the main empirical central limit theorems in \cite{vaartWellner2023}, as in Example~3.10.20.

Now it holds that
\begin{align*}
     &\sqrt{n}(\widehat{F}_1(t) - F_1(t)) \\&= \sqrt{n}\left(\int_0^t \widehat{S}(u-)\mathrm{d}\widehat{A}_1(u)  - \int_0^t {S}(u-)\mathrm{d}{A}_1(u)\right)
     \\&= \int_0^t \widehat{S}(u-)\mathrm{d}\sqrt{n}(\widehat{A}_1 - A_1)(u)  + \int_0^t \sqrt{n}(\widehat{S} - {S})(u-)\mathrm{d}{A}_1(u)
     \\&= \sqrt{n}(\widehat{A}_1 - A_1)(t)\widehat{S}(t) \!- \!\int_0^t \sqrt{n}(\widehat{A}_1 - A_1)(u)\mathrm{d}\widehat{S}(u) \! +\! \int_0^t \sqrt{n}(\widehat{S} - {S})(u-)\mathrm{d}{A}_1(u)
\end{align*}
 for all $t\in[0,K]$ by integration by parts, that is $$  \int_0^t f(v-) \;\mathrm{d}g(v) = (gf)(t) - (gf)(0-) - \int_0^t g(v) \;\mathrm{d}f(v)  $$ for $f\in BV_1[0,K], g\in D[0,K]$, where $BV_1[0,K]$ denote the set of all càdlàg functions $D[0,K]$ of total variation bounded by 1. 
As in Example~3.10.33 in \cite{vaartWellner2023}, the functional delta method yields
\begin{align*}
    \left( \sqrt{n}(\widehat{A}_1 - A_1) , \sqrt{n}(\widehat S - S)\right) \xrightarrow{d} \left( U_1, - S(.) \int_0^{.} \frac{S(v-)}{S(v)} \mathrm{d}(U_1 + U_2)(v)  \right)
\end{align*} as $n\to\infty$ on $D^2[0,K] $, where the integral is defined by integration by parts
since $U_1 + U_2$ is not of bounded variation.
Hence, we get \begin{align}\label{eq:slutsky}
\!\!\!    \left( \sqrt{n}(\widehat{A}_1 - A_1) , \sqrt{n}(\widehat S - S), \widehat{S} \right) \xrightarrow{d} \left( U_1, - S(.) \int_0^{.} \frac{S(v-)}{S(v)} \mathrm{d}(U_1 + U_2)(v), S  \right)
\end{align} as $n\to\infty$ on $D^2[0,K] \times BV_1[0,K] $
by Slutsky's lemma. 
Note that the map
\begin{align*}
\psi: D^2[0,K] \times BV_1[0,K] & \rightarrow  D[0,K], \\
      (\Tilde{A},\Tilde{B},\Tilde{C}) & \mapsto \Tilde{A}(.)\Tilde{C}(.) - \int_0^. \Tilde{A}\mathrm{d}\Tilde{C}  - \int_0^. \Tilde{B}(u-) \mathrm{d}{A}_1(u)
\end{align*}
is continuous on $D^2[0,K] \times \{S\}$.
Thus, an application of the continuous mapping theorem and changing the order of integration result in
\begin{align*}
    &\sqrt{n}(\widehat{F}_1 - F_1) \\&\xrightarrow{d}
    U_1(.)S(.) - \int_0^. U_1\mathrm{d}S  - \int_0^.  S(u-) \int_0^{u-} \frac{S(v-)}{S(v)} \mathrm{d}(U_1 + U_2)(v)  \mathrm{d}{A}_1(u)
    \\& =\int_0^. {S}(u-)\mathrm{d}U_1(u)  - \int_{[0,.)}\frac{S(v-)}{S(v)} \int_{(v,.]} S(u-)  \mathrm{d}{A}_1(u)\mathrm{d}(U_1 + U_2)(v)
    \\& =\int_0^. {S}(u-)\mathrm{d}U_1(u)  - \int_{[0,.)}\frac{S(v-)}{S(v)} (F_1(.) - F_1(v)) \mathrm{d}(U_1 + U_2)(v)
    \\& =\int_0^. \frac{S(u-)}{S(u)}\left({S}(u) - F_1(.) + F_1(u)\right)\mathrm{d}U_1(u)  + \int_0^. \frac{F_1(u) - F_1(.)}{1 - \Delta A(u)} \mathrm{d}U_2(u)
    \\& =\int_0^. \frac{1 - F_2(u) - F_1(.)}{1 - \Delta A(u)}\mathrm{d}U_1(u)  + \int_0^. \frac{F_1(u) - F_1(.)}{1 - \Delta A(u)} \mathrm{d}U_2(u)
\end{align*} as $n\to\infty$ on $D[0,K] $.

\begin{thm}[Corrected Theorem~4.1]
    As $n\to\infty$, we have on the càdlàg space $D[0,K]$
    \begin{align*}
        \sqrt{n}(\widehat{F}_1 - F_1) &\xrightarrow{d} U_{F_1} = \int_0^. \frac{1 - F_2(u) - F_1(.)}{1 - \Delta A(u)}\mathrm{d}U_1(u)  + \int_0^. \frac{F_1(u) - F_1(.)}{1 - \Delta A(u)} \mathrm{d}U_2(u),
    \end{align*}
    where $U_{F_1}$ is a zero-mean Gaussian process with covariance function
    \begin{align*}
        \sigma_{F_1}^2 : (s,t) \mapsto & \int_0^{s \wedge t} \frac{(1 - F_2(u) - F_1(s))(1 - F_2(u) - F_1(t))}{\bar{H}(u)}\frac{\mathrm{d}A_1(u)}{1 - \Delta A(u)}\\
    &+ \int_0^{s \wedge t} \frac{(F_1(u) - F_1(s))(F_1(u) - F_1(t))}{\bar{H}(u)}\frac{\mathrm{d}A_2(u)}{1 - \Delta A(u)}\\
    &+ \sum_{u \in D, u\leq s,t} \frac{S^2(u)}{\bar{H}(u)}\frac{\Delta A_1(u)\Delta A_2(u)}{(1 - \Delta A(u))^2}.
    \end{align*}
\end{thm}
The covariance function can be calculated analogously to Appendix~E.
Here, the last sum may be simplified to $ \sum_{u \in D, u\leq s,t} \frac{S(u-)}{G(u-)}\Delta A_1(u)\Delta A_2(u)$.



\section*{Funding}
Merle Munko gratefully acknowledges support from the \textit{Deutsche Forschungsgemeinschaft} (grant no. DI 2906/1-2 and GRK 2297 MathCoRe).

\end{document}